\documentclass[a4paper]{article}

\usepackage{
amsmath,
amsthm,
amscd,
amssymb,
wasysym
}
\usepackage{graphicx}
\usepackage{comment}
\usepackage{tikz}
\usetikzlibrary{decorations.text}
\usepackage{xspace}

\usepackage[backend=biber, style=numeric]{biblatex}
\usepackage{graphicx}

\usepackage{url}

\theoremstyle{plain}
\newtheorem{theorem}{Theorem}[section]
\newtheorem{lemma}[theorem]{Lemma}
\newtheorem{definition}[theorem]{Definition}
\newtheorem{corollary}[theorem]{Corollary}
\newtheorem{proposition}[theorem]{Proposition}
\newtheorem{conjecture}[theorem]{Conjecture}
\newtheorem{question}[theorem]{Question}
\newtheorem{remark}[theorem]{Remark}

\newtheorem{claim}[theorem]{Claim}

\newcommand*{\claimproofname}{Proof of claim}
\newenvironment{claimproof}[1][\claimproofname]{\begin{proof}[#1]}{\end{proof}}

\newtheoremstyle{derp}
{3pt}
{3pt}
{}
{}
{\upshape}
{:}
{.5em}
{}
\theoremstyle{derp}
\newtheorem{example}{Example}

\newcommand{\R}{\mathbb{R}}
\newcommand{\Q}{\mathbb{Q}}
\newcommand{\Z}{\mathbb{Z}}

\newcommand{\N}{\mathbb{N}}

\newcommand{\Fix}{\mathrm{Fix}}

\newcommand{\ID}{\mathrm{id}}

\newcommand{\vect}[1]{\overline{#1}}

\newcommand{\stab}{\mathrm{stab}}

\newcommand{\sign}{\mathrm{sgn}}

\newcommand\xqed[1]{%
  \leavevmode\unskip\penalty9999 \hbox{}\nobreak\hfill
  \quad\hbox{#1}}
\newcommand\qee{\xqed{$\fullmoon$}}

\newcommand{\id}{\mathrm{id}}
\newcommand{\pdim}{\mathrm{pdim}}

\newcommand{\BS}{\mathrm{BS}}

\title{Contractible subshifts}

\author{
Leo Poirier \and Ville Salo
}

\begin{document}
\maketitle

\begin{abstract}
We introduce the notion of a contractible subshift. This is a strengthening of the notion of strong irreducibility, where we require that the gluings are given by a block map. We show that a subshift is a retract of a full shift if and only if it is a contractible SFT with a fixed point. For many groups, including virtually polycyclic groups, metabelian Baumslag-Solitar groups and the lamplighter group, contractibility implies dense periodic points. We introduce a ``homotopy theory'' framework for working with this notion, and ``contractibility'' is in fact simply an analog of the usual contractibility in algebraic topology. We also explore the symbolic dynamical analogs of homotopy equivalence and strong contractibility of subshifts (corresponding to the notion of equiconnectedness in topology). Contractibility is implied by the map extension property of Meyerovitch, and among SFTs, it implies the finite extension property of Brice\~no, McGoff and Pavlov. We include thorough comparisons with these classes. We also encounter some new geometric notions, in particular a periodic variant of Gromov's asymptotic dimension of a group.
\end{abstract}

\section{Introduction}

If $G$ is a group, $\Sigma$ a finite set (alphabet) and $X \subset \Sigma^G$ is topologically closed, and closed under the shift maps induced by the regular action (given by the formula $(g \cdot x)_h = x_{g^{-1}h}$ where $g, h \in G, x \in X$) then $X$ is a \emph{subshift}. Subshifts are the main object of study in the field of symbolic dynamics.

An important notion in symbolic dynamics is \emph{gluing}. Suppose $X$ is a subshift, $x, y \in X$, and $A,B \subset G$ (if $A$ and $B$ are finite, we use the notation $A,B \Subset G$). Then a \emph{gluing} of $x,y$ on the respective areas $A, B$ is a point $z \in X$ such that $z|A = x|A$ and $z|B = y|B$. By a \emph{gluing notion} we understand a property of a subshift, which states the existence of gluings under some conditions.

The classical notion of topological transitivity of a $G$-system can be seen as a gluing notion. Recall that the definition of this is that if $U, V$ are nonempty open sets, then we can always find\footnote{Note that for $\Z$-systems, topological transitivity usually refers to the stronger variant, where recurrence is required in the positive direction. For general $G$, this has no natural counterpart, and it is common to use this weaker variant.} $x \in U$ and $g \in G$ such that $gx \in V$ (where the action of $n$ is written as $\sigma^n$). When we consider a $G$-subshift $X$, topological transitivity can be seen as a gluing property, as it is equivalent to the following statement: for any finite $A, B \Subset G$ and $x, y \in X$, there exists $g \in G$ such that $x$ and $gy$ can be glued on the areas $A$ and $gB$. 

Topological transitivity is too weak of an assumption for many symbolic dynamical purposes, especially on groups other than $\Z$. A much stronger gluing property is \emph{strong irreducibility} (see e.g.\ \cite{CeCo12}). It states that a gluing can be performed for any two points $x, y$ whenever the minimal distance between an element of $A$ and an element of $B$ is bounded from below by some constant (depending on the subshift). 

Strong irreducibility is perhaps the strongest possible ``pure'' gluing notion. Yet, it is not strong enough for many purposes. In particular, it is an open problem whether $\Z^d$-subshifts of finite type (see Section~\ref{sec:Pre} for the definition) have dense periodic points for $d \geq 3$, and it was recently shown that general strongly irreducible $\Z^d$-subshifts do \emph{not} always have dense periodic points for $d \geq 2$ \cite{Ho24}. 

Thus, it makes sense to try to strengthen gluing notions by adding other ingredients. A natural possibility is to require that instead of gluings simply existing, we have some control on the gluing. For example, Bowen's classical notion of specification for one-dimensional systems \cite{Bo71} includes (slightly rephrasing) the requirement that if the areas $A$, $B$ are periodic, then the gluing can also be made periodic. 

One situation where we want very fine control on the gluing is when constructing block maps to a subshift (for example, automorphisms). In this case, it is typically not useful to know that a gluing exists, but it should have some predictable form -- in fact, it should be produced by a block map itself. This is the approach taken in the present paper. 

\subsection{Contractibility}

In the present paper, we consider the notion of strong irreducibility with the requirement that the gluing is produced by a block map. This means that we should have a shift-commuting continuous map $h$ which is given two points $x, y$ from the subshift $X$, and another ``time'' configuration $t$ that designates areas where we glue from the first point and areas where we glue from the second point. Given this data, it produces a gluing that is valid sufficiently deep inside the given areas.

More precisely, letting $I$ denote the full shift $\{0,1\}^G$, we should have a \emph{block map} (meaning shift-commuting continuous map)
\[ h : I \times X \times X \to X \]
such that $h(t, x, y)$ copies the symbol from $x$ (resp.\ $y$) at $a \in G$ if $t$ contains a large patch of $0$s (resp.\ $1$s) near position $a$. Of course, by compactness of the Cantor space $I \times X \times X$, this is precisely equivalent to the requirements $h(\bar 0, x, y) = x$, $h(\bar 1, x, y) = y$, where $\bar a$ is the constant-$a$ configuration.

Note that $h(\bar 0, x, y) = x$, $h(\bar 1, x, y) = y$ is analogous to the usual definition of a homotopy in algebraic topology between two projection maps, except the interval is replaced by a full shift, and in addition to continuity we require shift-commutation.

\begin{definition}
\label{def:Homotopy}
If $f, g : X \to Y$ are morphisms (= block maps), a \emph{homotopy} between $f$ and $g$, denoted $h : f \cong g$, is a block map $h : I \times X \to Y$ satisfying $h(\bar 0, x) = f(x)$ and $h(\bar 1, x) = g(x)$. We then say $f$ and $g$ are \emph{homotopic}.
\end{definition}

We can then mimic (one of the standard definitions of) the notion of contractibility from algebraic topology.

\begin{definition}
A subshift $X \subset \Sigma^G$ is \emph{contractible} if there is a homotopy between the two projections $\pi_1, \pi_2 : X \times X \to X$.
\end{definition}

We show in Section~\ref{sec:ContractibilityInAT} that this corresponds to the standard notion of contractibility in algebraic topology.

Unwrapping the definition, we see that this means precisely what we wrote above: a subshift $X$ is contractible if there is a continuous, shift-commuting map $h : I \times X \times X \to X$ such that $h(\bar 0, x, y) = x$, and $h(\bar 1, x, y) = y$. We call such $h$ a \emph{contraction homotopy}, or a \emph{gluing morphism}.


\subsection{Results}

In this section, we highlight some of the results of the paper. First, in the most classical case of one-dimensional SFTs, we have the expected result:

\begin{theorem}
A one-dimensional subshift of finite type is topologically mixing if and only if it is contractible.
\end{theorem}

This is proved in Theorem~\ref{thm:MixingContractible}.

The following shows that contractible subshifts (even without the SFT assumption) avoid some pathologies of strongly irreducible subshifts.

\begin{theorem}
\label{thm:DPP}
Let $G$ be an infinite finitely-generated residually finite group with finite periodic asymptotic dimension. Then every contractible $G$-subshift has dense periodic points.
\end{theorem}

This is proved in Theorem~\ref{thm:ContractibleDPP}. Note that here ``periodic'' refers to having a finite index stabilizer.

Periodic asymptotic dimension is a new notion that we introduce. We show in Section~\ref{sec:Dimension} that strongly polycyclic groups (such as the group $\Z^d$), the lamplighter groups $\Z_p \wr \Z$, and the metabelian Baumslag-Solitar groups $BS(1, n)$, have finite periodic asymptotic dimension (and they are of course infinite finitely-generated residually finite groups).

In Section~\ref{sec:Retracts} we characterize retracts of full shifts in terms of contractibility, showing that contractible subshifts (of finite type, with a fixed point) arise naturally in the category of subshifts of finite type on any group (subshifts form a category in a natural way by taking morphisms to be (not necessarily onto) block maps \cite{SaTo15}).

\begin{theorem}
\label{thm:RetractResult}
On every group, a subshift $X$ is a retract of a full shift if and only if it has a fixed point, is of finite type, and is contractible.
\end{theorem}

This is proved in Theorem~\ref{thm:RetractCharacterization}.
An analogous result is more generally true when the full shift is replaced by an arbitrary contractible subshift $Y$, though the existence of a fixed point is replaced by the assumption that a map exists from $Y$ to $X$.

While our notion of contractibility is stated as a strengthening of strong irreducibility, its connection to the finite extension property (FEP) of \cite{BrMcPa18} seems much tighter. This is not particularly surprising, as the FEP was also introduced as a requirement on the target subshift of a block map. Recall that the finite extension property for an SFT means that if a pattern on a domain $D$ does not contain forbidden patterns, then the subpattern seen in cells sufficiently far from the boundary actually appears in a point of the SFT. See Section\ref{sec:Definitions} for the precise definition.

\begin{theorem}
On every group, every contractible SFT has the FEP.
\end{theorem}

This is proved in Theorem~\ref{thm:ContractibleSFThasFEP}. We show by examples in Section~\ref{sec:Comparisons} that on many groups, contractible SFTs are strictly contained in the class of FEP subshifts, which in turn are a strict subclass of factors of contractible SFTs.

The FEP class was introduced in \cite{BrMcPa18} as the codomain condition in a general result about factoring SFTs onto smaller entropy SFTs. In Theorem~\ref{thm:OntoContractible}, we prove the following weaker result, which applies on any group, and has a very simple proof.

\begin{theorem}
\label{thm:OntoContractibleResult}
Suppose $Y \subset A^G$ is a contractible SFT. Let $X \subset B^G$ be an arbitrary subshift. Then the following are equivalent:
\begin{itemize}
\item There is a block map from $f : X \to Y$ and a block map $g : X \to A^G$ such that $g(X) \supset Y$.
\item There is a factor map $f : X \to Y$.
\end{itemize}
\end{theorem}

Here, the first item lists the obvious necessary conditions for the existence of a factor map. It is not difficult to recover the special case of the result of \cite{BrMcPa18} for contractible SFTs from this theorem.

Contractibility and FEP are incomparable properties, and we do not know a natural common generalization, but there are many connections.
In particular in Theorem~\ref{thm:FEPDPP}, we prove the analog of Theorem~\ref{thm:DPP} for FEP subshifts, with essentially the same proof.

Another class of subshifts, partially motivated by similar concerns as our notion, called subshifts with the map extension property, was recently -- and indeed essentially simultaneously -- introduced by Meyerovitch in \cite{Me23}. At least on abelian groups, these properties are equivalent, except that the definition of Meyerovitch forces commutativity:


\begin{theorem}
\label{thm:ContractibleMEPCharacterization}
On $\Z^d$, contractible SFTs are precisely the subshifts satisfying the map extension property of \cite{Me23}.
\end{theorem}

This is proved in Theorem~\ref{thm:ZdFirstItemCharacterization}. We also give a characterization of this property on a general group in Theorem~\ref{thm:FirstItemCharacterization}, using an auxiliary notion we call the map existence property.

There is some overlap between \cite{Me23} and the present paper, partially due to their concurrency, and partially because ideas from each paper influenced the other, yet both papers were kept self-contained. More discussion of connections can be found in the paper of Meyerovitch, and we recommend reading his paper in parallel to the present one as it provides a different point of view to basically the same concept, and has different emphasis. 


There are many other things one can do with the symbolic dynamical homotopy notion in Definition~\ref{def:Homotopy}, in particular we study strong contractibility (which corresponds to the notion of equiconnectedness in topology), meaning that there is a homotopy between the two projections of $X \times X$ which fixes the diagonal. 

For this, we define a new technical notion of a group, called the \emph{patching property}, in Section~\ref{sec:Patching}. We show that it is implied by finite asymptotic dimension, and by subexponential growth. Thus it includes the groups listed above, but also for example nonabelian free groups, and the Grigorchuk group. We do not have any example of a group that does not have the patching property.

\begin{theorem}
Suppose the group $G$ has the patching property. Then all strongly contractible subshifts are contractible SFTs. If the group is $\Z^d$, or the subshift has a fixed point, the converse holds.
\end{theorem}

The first direction is proved in Theorem~\ref{thm:EqContContSFT}. The second is proved in Theorem~\ref{thm:ZdFirstItemCharacterization} ($\Z^d$) and Theorem~\ref{thm:ContractibleFPEqConn} (subshifts with fixed points).

We also take a look at homotopy equivalence, and show among other things the following:

\begin{theorem}
Two one-dimensional transitive SFTs are homotopy equivalent if and only if they map into each other.
\end{theorem}

\begin{theorem}
A subshift of finite type is homotopy-equivalent to a point if and only if it is contractible and has a fixed-point.
\end{theorem}

The first is proved in Theorem~\ref{thm:HomotopyEquivalent}, and the latter in Theorem~\ref{thm:HEToPoint}. (Note that this is the same condition that characterizes the retracts of full shifts.)

In the important special case $G = \Z^d$, we summarize the equivalence results obtained above, first in the case of SFTs with a fixed point, and then in general:

\begin{theorem}
Let $G = \Z^d$, and let $X$ be a subshift of finite type with a fixed point. Then TFAE:
\begin{itemize}
\item $X$ is contractible;
\item $X$ is a retract of a full shift;
\item $X$ is a retract of a contractible subshift of finite type;
\item $X$ is strongly contractible;
\item $X$ satisfies the map extension property of \cite{Me23};
\item $X$ is homotopy-equivalent to a point.
\end{itemize}
\end{theorem}


\begin{theorem}
Let $G = \Z^d$, and let $X$ be a subshift. Then TFAE:
\begin{itemize}
\item $X$ is a contractible SFT;
\item $X$ is a retract of a contractible SFT;
\item $X$ is strongly contractible;
\item $X$ satisfies the map extension property of \cite{Me23}.
\end{itemize}
\end{theorem}

Besides specific results, we give many examples of contractible subshifts and general constructions throughout the paper.

\section{Preliminaries}
\label{sec:Pre}

\subsection{Contractibility in algebraic topology}
\label{sec:ContractibilityInAT}

Our definition of contractibility indeed corresponds to the usual notion in the setting of topological spaces. The most common definition of contractibility in algebraic topology is likely that a topological space is contractible if there is a homotopy between the identity map on $X$, and the map from $X$ to one of its points. We show below that in the topological setting, this is equivalent to the definition we use in this paper.

\begin{lemma}
\label{lem:TwoHomotopyDefinitions}
Let $X$ be a topogical space. Then the following are equivalent:
\begin{itemize}
\item There is a homotopy between the identity map on $X$, and the constant map from $X$ to one of its points.
\item The two projections from $X \times X$ to $X$ are homotopic.
\end{itemize}
\end{lemma}

\begin{proof}
Suppose there is a homotopy between the identity map on $X$, and the map from $X$ to one of its points, say $h : \id_X \cong (x \mapsto x_0)$ for some $x_0 \in X$. Then 
\[ h'(t, x, y) = \left\{
\begin{array}{ll}
h(2t, x) & \mbox{ if } t \leq 1/2 \\
h(2-2t, y) & \mbox{ if } t \geq 1/2 
\end{array}\right. \]
is a homotopy between the two projections. Conversely if $h'$ is a homotopy between the two projections, then
\[ h(t, x) = h'(t, x, x_0) \]
gives a homotopy between the identity map on $X$, and the map from $X$ to $\{x_0\}$.
\end{proof}

In our symbolic dynamical category, these notions are not equivalent. The basic issue being that a morphism cannot talk about any individual point unless there is a fixed point. We will say a subshift is \emph{fixed-point contractible} if the analog of the more common definition of contractibility holds, i.e.\ the identity map is homotopic to the map to some one-point subsystem $\{\bar a\}$. We will see that contractibility and fixed-point contractibility are equivalent concepts for subshifts of finite type containing at least one fixed point (Corollary~\ref{cor:FixedPointContractible}). We find the second definition more useful in our symbolic dynamical category, so we use that as the main definition. 

\subsection{General and symbolic dynamical definitions}
\label{sec:Definitions}

If $f : A \to B$ is a function and $C \subset A$, then $f|C$ is the domain restriction of $f$ to $C$. By $A \subset B$ we mean $A$ is a not necessarily proper subset of $B$. By $\Z_k$ we denote the finite ring or additive group $\Z/k\Z$. We write $A \subset B \iff A \subseteq B$.

Throughout, $G$ denotes a group. The identity element for abstract groups is $1_G$. Elements of a group are usually called $a, b, c, e$ ($f$, $g$ and $h$ will be morphisms, and $d$ is a distance). We often refer to elements of a group as \emph{cells}, a terminology coming from cellular automata.

To a group $G$ with symmetric generating set $S$, we associate Cayley graph of $G$ as having edges $(a, ab)$ for $b \in S$ (i.e.\ the right Cayley graph), which we sometimes take to have edge-label $b$. In this point of view, we often call elements of the group nodes.

We mostly consider infinite generated groups with a preferred generating set $S$. We usually (but not always) let $S$ be a finite symmetric ($S = S^{-1}$) generating set. We write $d(a, b)$ for the left-invariant word metric of $G$ with respect to this generating set. The \emph{ball of radius} $R$ is $B_R = \{g \in G \;|\; d(g, e) \leq R\}$. It is a metric ball in the right Cayley graph, with the path metric.

When we consider a subset $H$ of a group $G$, we think of it as a metric space with the induced metric from $G$. Note that this may not be the metric induced by the path metric on the induced subgraph $H$ with respect to any generating set of $G$. We say a subset $H \subset G$ is \emph{$r$-connected} if, defining graph $(H, E)$ where we have vertices $H$, and edges $(h, h') \in E$ when $d(h, h') \leq r$ (in the group $G$), the graph $(H, E)$ is connected. More generally, we define the \emph{$r$-components} as the connected components of this graph.

When $G \curvearrowright X$ is an action and $H \leq G$, write $\Fix(H)$ for the set of fixed-points of the action of $H$. We also call these the \emph{$H$-periodic points}.

The \emph{full shift} is the set $A^G$ (functions from $G$ to $A$), for a finite set $A$ called the \emph{alphabet}, whose elements are called \emph{symbols} or sometimes \emph{colors}. Its elements $x \in A^G$ are called \emph{configurations} (or sometimes \emph{points} to emphasize the dynamical systems aspect) and we write $x_a$ instead of $x(a)$, but we use the same restriction notation as for functions, and on occasion the notation $x^{-1}(a)$ for the set of cells in $G$ where $x$ contains the symbol $a$. A full shift $A^G$ is \emph{non-trivial} if $|A| \geq 2$. 

The full shift is a topological dynamical system when $A^G$ is given the product topology (making it a Cantor set if $|A|\geq 2$). $G$ acts by $ax_b = x_{a^{-1}b}$ for $x \in A^G$ and $a, b \in G$. Under this shift convention, if we see a configuration as a vertex-coloring of the $S$-labeled right Cayley graph of the group (where $S$ is a generating set), then the shift map by $a$ translates the identity $1_G$ to $a$, and translates everything else by the unique graph automorphism of the Cayley graph.

If $K \leq G$ is a subgroup, a configuration $x \in A^G$ is \emph{$K$-periodic} if its stabilizer contains $K$, in other words if $ax = x$ for all $a \in K$. A point is \emph{periodic} is it is $K$-periodic for a finite-index subgroup $K \leq G$.\footnote{Such a point is often called ``totally periodic'' in symbolic dynamics, but since having just some nontrivial period is rarely of interest in the case of subshifts on general groups, the word ``periodic'' usually refers to the finite index case in the present paper, when used in compound terms.} Being periodic is equivalent to having a finite orbit under the $G$-action.

If $D \subset G$, and $p \in A^D$, then $p$ is a \emph{pattern}. A pattern is \emph{finite} if its domain $D$ is finite, which we denote by $D \Subset G$. We translate patterns as follows: if $p \in A^D$ and $a \in G$, then $ap \in A^{aD}$ and $ap_b = p_{a^{-1}b}$.
Note that configurations are patterns with a full domain, and their translation agrees with the translation of patterns.
If $x \in A^G$ then our syntactic convention is that translation happens \emph{before} restriction, so $ax|N = (ax)|N$. We note the useful formula $a^{-1}(x|aN) = a^{-1}x|N$.

If $p, q$ are patterns, then $q \sqsubset p$ where $q \in A^D$, $p \in A^E$, means $p|aD = aq$ for some $a \in G$ such that $aD \subset E$. We say $q$ \emph{appears} in $p$ or $p$ \emph{contains} $q$. We extend this naturally to the case where $p \in A^E$ for possibly infinite $E \subset G$ (e.g.\ a $p$ a configuration in some subshift).

We extend this to the case where $q$ is a single symbol $s$, which we can identify with the pattern $1_G \mapsto s$. We write $q \subset p$ if $q$ is a subpattern of $p$ as a positioned pattern, meaning $p|D = q$ where $q$ has domain $D$. If $X \subset A^D$ is a set of patterns (possibly with infinite $D$), we say $p$ \emph{appears} in $X$, or $X$ \emph{contains} $p$, if $p$ appears in $x$ for some $x \in X$.

The \emph{cylinder} corresponding to a finite pattern $p$ is $[p] = \{x \in A^G \;|\; \forall a \in D: x_a = p_a\}$. The topology of $A^G$ has a basis of clopen sets. The clopen sets are precisely the finite unions of \emph{cylinders}.

A \emph{subshift} is a topologically closed and $G$-closed subset $X$ of a full shift. Equivalently, given a countable (possibly finite) set $\mathcal{F} := \{p_1, p_2, \ldots \}$ of finite patterns, it is a subset of $A^G$ defined as the set of points $x \in A^G$ whose \emph{orbit} $Gx$ does not intersect any cylinder $[p_i]$ for $p_i \in \mathcal{F}$. The corresponding $p_i$ are called \emph{forbidden patterns} of $X$. This is equivalent to saying that $X$ is the set of points in $A^G$ which do not contain any of the forbidden patterns in $\mathcal{F}$.

If we can pick a finite set of forbidden patterns, then $X$ is a \emph{subshift of finite type} or \emph{SFT} for shift. We may always suppose the forbidden patterns of an SFT have the same finite domain $D$, which we call a \emph{window} for the SFT. If only the window $W$ is specified, but not the patterns, we implicitly mean that the forbidden patterns are precisely the patterns of shape $W$ which do not appear in configurations. Note that by our shift convention, forbidden patterns are checked in the right Cayley graph, i.e.\ for each $g$, we check a condition on the pattern $x|gD$.

If $X \subset Y$ are subshifts, we say $X$ is a \emph{relative SFT} of $Y$ if $X$ is the intersection of $Y$ with an SFT, equivalently it admits a defining set of forbidden patterns consisting of those of $Y$, and finitely many additional ones.

When we have fixed a finite set of forbidden patterns $\mathcal{F} \subset A^D$ for $D \Subset G$, defining an SFT $X \subset A^G$, we say a pattern $p \in A^E$ is \emph{locally valid} in $X$ if it does not contain any of the patterns in $\mathcal{F}$. 
We say it is \emph{globally valid} if it can be extended to a valid configuration of $X$. This means precisely that $p$ appears in $X$, i.e.\ $p \in X|D$.

A subshift $X \subset \Sigma^G$ is \emph{strongly irreducible} if there exists $R$ such that for all $A, B \subset G$ such that $d(A, B) = \min_{a \in A, b \in B} d(a, b) \geq R$, and for all $x, y \in X$, there exists $z \in X$ such that $z|A = x|A, z|B = y|B$. We say $X$ is \emph{(topologically) mixing} if for all $A, B \Subset G$ (finite subsets!) there exists $R$ such that whenever $x, y \in X$, and $a \in G$ is such that $d(A, aB) \geq R$, there exists $z \in X$ such that $z|A = x|A, z|aB = y|aB$. It is easy to show that contractibility implies strong irreducibility. Lemma~\ref{lem:ContractibleSI} shows a stronger version of this result.

The \emph{finite extension property} or \emph{FEP} \cite{BrMcPa18} for an SFT $X$ means that for some SFT window $W \Subset G$, there exists $1_G \in N \Subset G$ such that if a pattern $p \in \Sigma^D$ admits a locally valid extension to $DN$, then $p$ is globally valid. We call $N$ the \emph{FEP window}. Note that an FEP subshift is in particular SFT with $W$.

Especially when we have homotopy in mind, we denote by $I_n$ the full shift on symbols $\{0,1,\ldots,n-1\}$. The binary full shift $I = I_2$ plays the role of the interval, and we often refer to its points as \emph{time parameters}. For a symbol $i$, write $\bar i$ or sometimes $i^G$ for the unique point in $\{i\}^G$. For a function $p : D \to A$ we also use the usual notation $p \equiv i$ for $\forall a \in D: p_a = i$, where $i \in A$. So $\bar i$ is the unique configuration $x$ with $x \equiv i$. Such configurations are often called \emph{unary}. We also use such terminology and notation for patterns, writing $i^D$ for the unique element for $\{i\}^D$.

We work in the category of subshifts on $G$, and a \emph{morphism} or \emph{block map} is a shift-commuting continuous function between two subshifts. Isomorphisms, embeddings and factor relations are naturally interpreted with these morphisms, as bijective, injective and surjective morphisms, respectively. The image of a surjective morphism is also called a \emph{factor}, and a preimage in a surjective morphism is a \emph{cover}. Isomorphisms are often called \emph{conjugacies}. Almost all of our definitions and notions are preserved under conjugacy (with some obvious exceptions like safe symbols and the Z0 property).

In a few occasions, it is convenient to consider more general systems than subshifts. We can extend the notion of conjugacy to compact metrizable spaces with continuous group actions by a countable group $G$. It is well-known that such a $G$-system $X$ is conjugate to a subshift if and only if $X$ is homeomorphic to a closed subset of Cantor space and the action is \emph{expansive} meaning \[ \exists \epsilon > 0: \forall x, y \in X: x \neq y \implies \exists g \in G: d(g(x), g(y)) > \epsilon. \]

Being conjugate to an SFT also has such a ``dynamical'' interpretation for general systems $(G, X)$: an \emph{$(\epsilon, S)$-pseudo-orbit}, where $\epsilon > 0$ and $S \Subset G$, is a map $p : G \to X$ such that $d(gx, sgx) < \epsilon$ for all $x \in X, s \in S, g \in G$. We say $x \in X$ \emph{$\delta$-shadows} $p$ if $d(gx, p(g)) < \delta$ for all $g \in G$. We say $X$ has the \emph{shadowing property} (a.k.a.\ pseudo-orbit tracing property) if for all $\delta > 0$, there exist $\epsilon > 0$ and $S \Subset G$ such that all $(\epsilon, S)$-pseudo-orbits are $\delta$-shadowed by some point of $X$. Now, a general system is conjugate to an SFT if and only if it is conjugate to a subshift, and has the shadowing property.

A morphism $f : X \to Y$ between two subshifts $X \subset A^G, Y \subset B^G$ has a \emph{neighborhood} $N \Subset G$ and a \emph{local rule} $f_{\mathrm{loc}} : A^N \to B$, such that $f(x)_a = f_{\mathrm{loc}}(a^{-1}(x|aN))$ (or equivalently $f(x)_a = f_{\mathrm{loc}}(a^{-1}x|N))$). We usually deal with only finitely-generated groups with a preferred generating set $S$, and then we also say $r$ is a \emph{radius} for $f$ it the radius-$r$ ball $B_r = S^{\leq r} = \{g \;|\; d(1_G, g) \leq r\}$ is a neighborhood for it.


If $X \subset A^H$ is a subshift and $\pi : G \to H$ is a group epimorphism, then the \emph{pullback} (of $X$ with respect to $\pi$) is the subshift $Y \subset A^G$ with
\[ y \in Y \iff \exists x \in X: \forall a \in G: y_a = x_{\pi(a)}. \]

If $H \leq G$, then corresponding the \emph{free extension} (a.k.a.\ \emph{induction}) of a subshift $X \subset A^H$ is the subshift $Y \subset A^G$ with the same forbidden patterns. We write it as $Y = X^{G/H}$ Equivalently, configurations of $Y$ consist of independently chosen configurations of $X$ on different left cosets of $H$.

If $X \subset A^G$ is a subshift, the \emph{SFT approximation} of $X$ with \emph{window} $W \Subset G$ is the SFT defined by forbidden patterns $\{ p \in A^W \;|\; [p] \cap X = \emptyset \}$. SFT approximations are ordered by \emph{goodness}: if $X',X''$ are SFT approximations of $X$ defined with windows $W', W''$ respectively, we say $X'$ is a \emph{better} approximation than $X''$ if $W'' \subset W'$.

If $A \times B$ is a (Cartesian) product alphabet, we think of subshifts $X \subset (A \times B)^G$ also as subshifts of $A^G \times B^G$ with the diagonal action. We often refer to the different projections of subshifts (and configurations) as \emph{tracks}. 

A \emph{vertex shift} is a subshift of $A^\Z$ defined by forbidden patterns with domain $\{0,1\}$. In the one-dimensional case (i.e.\ the group $\Z$) we also call patterns whose domain is an interval $\{0, 1, \ldots, k-1\}$ \emph{words}, and the \emph{length} of a word $w \in A^k = A^{\{0, 1,\ldots, k-1\}}$ is $|w| = k$. We identify them with the free monoid on the alphabet. We can \emph{concatenate} two words $u \in A^k, v \in A^m$ in an obvious way to a word $uv \in A^{k+m}$. For $A$ an alphabet, write $A^*$ for the set of all words over this alphabet, including the unique word of length $0$.

\subsection{A few words on groups}

We mostly assume the reader is familiar with group theory, but we recall some definitions and a basic lemma.

A \emph{polycyclic group} is $G$ such that there is a sequence of subgroups $1 = G_0 < G_1 < \ldots < G_n = G$ such that each $G_i$ is normal in $G_{i+1}$ and each quotient $G_{i+1}/G_i$ is cyclic. If the quotients are isomorphic to $\Z$, the group is \emph{strongly polycyclic}.

A group is \emph{just-infinite} if it is infinite and has no non-trivial normal subgroup of infinite index. Just-infinite groups are an important and interesting class of groups, see \cite{BaRoZo03}. Nevertheless, arguably most groups are not just-infinite. The group $\Z$ is the only infinite strongly polycyclic just-infinite group. 

We write $1 \rightarrow K \rightarrow G \rightarrow H \rightarrow 1$ (an \emph{exact sequence}) and say $G$ is a \emph{group extension} (of $H$ by $K$) if there is a surjective group homomorphism $\pi : G \to H$, such that $K$ is the kernel of this homomorphism. If $\pi$ has a section $\gamma : H \to G$, meaning $\gamma$ is a homomorphism and $\pi\circ\gamma = \id_H$, then we say $\pi$ (or the group extension) \emph{splits}, and then $G$ is a semidirect product $K \rtimes H$.

We sometimes need long or infinite geodesics. We recall the simple proof that these exist universally.

\begin{lemma}
\label{lem:Geodesics}
Let $G$ be infinite and finitely-generated by a symmetric set $S \Subset G$. Then there exists $p : \Z \to G$ such that $d(p(i),  p(j)) = |j - i|$ for all $i < j$.
\end{lemma}

\begin{proof}
By K\H{o}nig's lemma we have $p : \N \to G$ with this property, by considering optimal paths to the sphere $S^r$ for increasing $r$. Then for each $n$, reparametrize the path to $p_n(j) : [-n, \infty) \to G$ by $p_n(j) = p(n)^{-1} p(j + n)$, extend each $p_n$ arbitrarily to domain $\Z$, and take any limit point of the $p_n$ in the product topology.
\end{proof}

\subsection{Two lemmas about SFTs}

\begin{lemma}
\label{lem:SFTCharacterization}
The following are equivalent for a subshift $Y \subset B^G$:
\begin{enumerate}
\item $Y$ is SFT.
\item For any local rule $f_{\mathrm{loc}} : A^N \to B$ defining a valid map $f : X \to Y$ for a subshift $X \subset A^G$, for all sufficiently good SFT approximations $X'$ of $X$, $f_{\mathrm{loc}}$ also gives a well-defined map $\tilde f : X' \to Y$. 
\item Whenever $f : B^G \to B^G$ defined by $f_{\mathrm{loc}} : B^N \to B$ restricts to the identity map $\id : Y \to Y$, for all good enough SFT approximations $Y'$ of $Y$, we have $f(Y') \subset Y$. 
\end{enumerate}
\end{lemma}

\begin{proof}
For (1) to (2), 
let $W$ be a window for $Y$, and let $X'$ be the SFT approximation of $X$ with window $WN$ or greater (with respect to inclusion).

Now when applied to configurations of $X'$, $\tilde f$ sees the same $WN$-patterns as when applied to configurations of $X$. Thus, it produces the same $W$-patterns in the image, in particular if $f(X) \subset Y$ then it cannot produce forbidden patterns of $Y$. 

The implication (2) $\implies$ (3) is trivial.

For (3) to (1), let $f$ have neighborhood $N$. Let $Y'$ be a good enough SFT approximation that $f(Y') \subset Y$, and whose defining window contains $N$. Then $f$ is the identity map on $Y'$, so $Y = Y'$ is SFT.
\end{proof}

The following is a simple compactness argument.

\begin{lemma}
\label{lem:LocalGlobal}
Let $X \subset A^G$ be a subshift of finite type with a fixed set of forbidden patterns $P \subset A^K$. Then for all $T_1 \Subset G$, there exists $T_2 \Subset G$ such that when $p \in A^{T_2}$ contains no pattern from $P$, then $p|T_1$ appears in a configuration of $X$.
\end{lemma}


\begin{remark}
We note that the previous lemma is only an existential result, in that for the groups $G = \Z^d$, $d \geq 2$ (and more generally a large class of one-ended groups), there is no computable function that produces $T_2$ from $T_1$ and the forbidden patterns defining $X$. This is because we can encode seeded Turing machine computation into $T_1$-patterns, and the size of the corresponding $T_2$ roughly corresponds to the busy beaver function from computability theory.
\end{remark}

\subsection{The natural relaxation of a contraction homotopy}

A very common modification we need to perform on a contraction homotopy is to extend its local rule, so that it can also be applied in invalid contexts. We do this in a specific way.

\begin{definition}
\label{def:NaturalExtension}
Let $h : I \times X \times X \to X$ be a contraction homotopy, where $X \subset A^G$. Then a \emph{natural relaxation} is any $\tilde h : I \times A^G \times A^G \to A^G$ with the same neighborhood as $h$ such that $\tilde h(\bar 0, x, y) = x$, $\tilde h(\bar 1, x, y) = y$ for all $x, y \in A^G$ and $\tilde h(t, x, y) = h(t, x, y)$ whenever $x, y \in X$.
\end{definition}

A natural relaxation always exists: simply use the same local rule when you can, extend it to copy from the second input when the first input is all $0$, from the third when it is all $1$, and finally extend arbitrarily to other patterns.

\section{First examples of contractible SFTs}

An important definition, and the simplest possible example of a homotopy, is the following.

\begin{definition}
\label{def:NaiveHomotopy}
Let $G$ be any group, $X =\Sigma^G$, and $f, g : X \to X$ arbitrary block maps. For $a \in G, t \in I, x \in X$, define
\[ h(t, x)_a = \left\{\begin{array}{ll}
f(x)_a & \mbox{ if } t_a = 0, \\
g(x)_a & \mbox{ if } t_a = 1. \\
\end{array}\right. \]
Then $h$ is the \emph{naive homotopy} between $f$ and $g$.
\end{definition}

More generally, we refer to any homotopy defined by this formula a naive homotopy.

An important source of strongly irreducible subshifts are safe symbols. Recall that a \emph{safe symbol} for a subshift $X \subset A^G$ is a symbol $0 \in A$ such that for all $a \in G$ and $x \in X$ we have $y \in X$, where $y_b = x_b$ for $b \neq a$, and $y_a = 0$. Note that by compactness, any subset of symbols may be changed to $0$ in any valid configuration, and this always results in a valid configuration. 

In particular, every \emph{hereditary subshift}, meaning one where the alphabet $A$ is totally ordered and any symbol can be replaced by a smaller one, in any valid configuration, without introducing a forbidden pattern, has safe symbol $\min A$. A concrete example (which is also SFT so the previous proposition applies) is the golden mean shift in one dimension, i.e.\ the set $X \subset \{0,1\}^\Z$ of configurations $x \in \{0,1\}^\Z$ such that $11 \not\sqsubset x$.

\begin{proposition}
Every subshift of finite type which has a safe symbol is contractible.
\end{proposition}

\begin{proof}
Let $F\Subset G$ be a window for the SFT $X \subset A^G$ and suppose $0 \in A$ is a safe symbol. Let $M = F^{-1}$ and define $h : I \times X \times X \to X$ by
\[ h(t, x, y)_a = \begin{cases}
0 & \mbox{if } t|{aM} \notin \{0^{aM}, 1^{aM}\}, \\
x_a & \mbox{if } t|{aM} = 0^{aM}, \\
y_a & \mbox{if } t|{aM} = 1^{aM}.
\end{cases}
\]


Given inputs $x$ and $y$ in $X$, the only possibility for $h(t,x,y)|F$ to be forbidden is that both of the two latter cases of the definition are used when computing it. Namely, otherwise $h(t, x, y)|F$ is obtained from $x|F$ or $y|F$ by inserting safe symbols and thus cannot be forbidden.

So for a forbidden pattern to appear, we must have two cells $a,b \in F$ such that $t|aM = 0^{aM}$ and $t|bM = 1^{bM}$. But since $M = F^{-1}$, both $aM$ and $bM$ contain the cell $1_G$, which gives us a contradiction since it cannot be mapped by $t$ to $0$ and $1$ at the same time.
\end{proof}

Without the assumption of finite type, safe symbols do not imply contractibility. In Definition~\ref{def:Pk}, we define the subshift $P_2$, which is a union of two full shifts sharing a unary point $\bar 1$. It has safe symbol $1$, but is not contractible. See also the discussion after Corollary \ref{cor:FixedPointContractible}.

One can produce many examples without safe symbols by conjugating a subshift with a safe symbol to an isomorphic one without a safe symbol.

The following proposition gives an example of a contractible subshift of finite type that is not produced this way (though we do not prove that it cannot be conjugated to have a safe symbol). This example is from \cite{BuSt94}, and it was the first example of a strongly irreducible subshift of finite type with multiple measures of maximal entropy. We thus state it as a proposition (although it gives examples of contractible SFTs on any group).

\begin{proposition}
On $\Z^d$, $d \geq 2$, there exist contractible subshifts of finite type with no safe symbol and with multiple measures of maximal entropy.
\end{proposition}

Thus, contractibility (even with the additional assumption of SFTness) does not remove the a priori ``pathological'' possibility of multiple measures of maximal entropy. On the other hand, the example from \cite{BuSt94} is extremely natural, so including it is difficult to see it as a shortcoming of our definition.

\begin{proof}
Let $G$ be a group with symmetric generating set $S \Subset G \setminus \{1_G\}$. Consider the $G$-SFT $X$ over alphabet $\Sigma = \{-M, \ldots, -1, 1, \ldots, M\}$ where for each $s \in S$, we forbid $p \in \Sigma^{\{1_G, s\}}$ when $\prod_i p_i \leq -2$. In other words, adjacent symbols cannot have different signs, unless both are $\pm 1$.

For $G = \Z^d$, this is Example~1.5 of \cite{BuSt94}. They show that if $M$ is large enough and $d \geq 2$, there are multiple measures of maximal entropy for $X$.

By definition, $X$ is a subshift of finite type. It has no safe symbol (for $M \geq 2$), as the symbol $s \in \Sigma$ cannot be positioned next to $-\sign(s) 2$. 

The subshift $X$ is contractible. Namely, given $x,y \in X$ and $t \in I$, we define
\[ h(t, x, y)_a = \begin{cases}
x_a & \mbox{if } t|aS = 0^{aS} \\
y_a & \mbox{if } t|aS = 1^{aS} \\
\sign(x_{ab}) & \mbox{if } t|aS \notin \{0^{aS}, 1^{aS}\} \wedge \exists b \in S: t|abS = 0^{abS} \\
\sign(y_{ab}) & \mbox{if } t|aS \notin \{0^{aS}, 1^{aS}\} \wedge \exists b \in S: t|abS = 1^{abS} \\
1 & \mbox{otherwise.}
\end{cases} \]

This is well-defined, since if $\exists b \in S: t|abS = 0^{abS}$ then $t_a = 0$, so the third and fourth cases do not collide.
We claim that $h(t, x, y) = z$ can not contain a forbidden pattern. Otherwise, by shifting we may suppose $z_{1_G} z_a \leq -2$ for $a \in S$. Then certainly for $T = \{1_G, a\}$, $z|T \neq x|T$ and $z|T \neq y|T$, so $t|TS$ is non-unary.

If neither $t|S$ nor $t|aS$ is monochromatic, then both entries of $z|T$ have absolute value $1$. We conclude that exactly one of $t|S, t|aS$ is unary. By symmetry we may assume $t|aS$ is unary, and by symmetry we may assume $t|aS = 0^{aS}$. Then $z_a = x_a$ and $z_{1_G} = \sign(x_a)$, so in fact $z_{1_G} z_a = |z_a| \geq 1$. This contradicts the assumption that $T$ contains a forbidden pattern.

This concludes the proof that $h$ is a contraction homotopy for $X$.
\end{proof}

The examples above all have fixed points. There are also contractible subshifts without fixed points. The following example also appears in \cite{Me23}. It is often referred to as the subshift of graph colorings, but in this paper colors refer to symbols in a general configuration, so we avoid this name.


\begin{proposition}
Let $G$ be a group and let $S \not\ni 1_G$ be a finite symmetric subset of $G$. Let $|A| \geq |S|+1$, and let $X \subset A^G$ be the subshift containing configurations $x\in A^G$ such that $x_a \neq x_{ab}$ for all $b \in S$. Then $X$ is a contractible SFT.
\end{proposition}

\begin{proof}
%
%
%
%
We define a contraction homotopy $h : I \times X \times X \to X$. Let $f : (A \cup \{\#\})^S \to A$ be any function such that $f(p) \not\sqsubset p$ (in the sense that the symbol $f(p)$ does not appear in $p$) for all $p \in (A \cup \{\#\})^S$. Let $A = \{1, 2, \ldots, k\}$.

Define $h_0 : I \times X \times X \to (A \cup \{\#\})^G$ by
\[ h_0(t, x, y)_ a = \begin{cases}
x_a & \mbox{ if } t_a = 0, \\
y_a & \mbox{ if } t|aS = 1^{aS}, \\
\# & \mbox{ otherwise}.
\end{cases} \]

Then for each $i \in A$ define $h_i : I \times X \times X \to (A \cup \{\#\})^G$ by
\[ h_i(t, x, y)_a = \begin{cases}
h_{i-1}(t, x, y)_a & \mbox{ if } h_{i-1}(t, x, y)_a \neq \#, \\
f(a^{-1}(h_{i-1}(t, x, y)|aS)) & \mbox{ if } h_{i-1}(t, x, y)_a = \# \wedge x_a = i, \\
\# & \mbox{ otherwise.}
\end{cases} \]
%

Note that each $h_i$ will preserve all non-$\#$ symbols picked in previous stages. The image of $h = h_k$ cannot contain symbol $\#$: if $x_a = i$, then the definition of $h_i(t, x, y)_a$ will set it to a non-$\#$ symbol unless it was already set in $h_{i-1}(t, x, y)$. The image of $h$ is in fact contained in $X$, as clearly no map $h_i$ will produce a forbidden pattern: $h_0$ copies all symbols from $x$ or $y$, and never picks the value of two adjacent cells from different points; other $h_i$ pick values by applying $f$ to the pattern in the $S$-adjacent cells (and $f$ picks some value that does not appear in the neighborhood), and $f$ is never applied in consecutive cells in the same step.

Also, $h(\bar 0, x, y) = x$ and $h(\bar 1, x, y) = y$ are already clear from the definition of $h_0$.
\end{proof}


More examples of contractible subshifts are given throughout the paper. In particular in Proposition~\ref{prop:PmSFTApproContractible}, in Lemma~\ref{lem:ContractibleDimensionSFT}, and several examples in Section~\ref{sec:Comparisons}.



\section{Characterization of retracts of full shifts}
\label{sec:Retracts}

One natural situation where contractible subshifts appear is the characterization of retracts of full shifts. 

\begin{definition}
Let $X, Y$ be two subshifts. We say $X$ is a \emph{retract} of $Y$ if there is a morphisms $f : X \to Y$ and a factor map $r : Y \to X$ such that $r \circ f = \id_X$.
\end{definition}



In this situation we say $f$ is a \emph{section} of $r$, and that $r$ is a \emph{retraction} of $f$. In this case one also says that $f$ is \emph{split monic} \cite{SaTo15}, though we do not use this terminology here. Note that any $f$ admitting a retraction needs to be injective, thus a conjugacy to its image, and it is easy to show that once $f$ is indeed injective, being split monic only depends on the image $f(X)$. In particular, when classifying retracts of a subshift $Y$, it suffices to classify retracts of $Y$ which are literally subshifts of $Y$, and $r|_X = f = \id_X$.

\begin{lemma}
For every subgroup $H \leq G$, if $Y$ has an $H$-periodic point so do its retracts.
\end{lemma}

\begin{proof}
The retraction map $r : Y \to X$ commutes with $G$ so $h y = y \implies r(h y) = h r(y)$. Thus if $H$ fixes $y$, it also fixes $r(y) \in X$.
\end{proof}

\begin{lemma}
A retract of a subshift of finite type is of finite type.
\end{lemma}

\begin{proof}
Suppose $Y$ is a subshift of finite type and suppose $f : X \to Y$ and $r : Y \to X$ with $r \circ f = \id_X$. Pick a local rule $g_{\mathrm{loc}}$ for $r \circ f$ (for example by composing the local rules of $f$ and $r$ in the obvious way). Since $Y$ is of finite type, by Lemma~\ref{lem:SFTCharacterization} for any good enough SFT approximation $X'$ of $X$, there is an extension $\tilde f : X' \to Y$ of $f$ with the same local rule.

Then $r \circ \tilde f : X' \to X$ is an extension of $\id_X$, and for good enough approximation uses the same local rule $g_{\mathrm{loc}}$. By Lemma~\ref{lem:SFTCharacterization}, $X$ is of finite type.
\end{proof}

\begin{lemma}
\label{lem:RetractContractible}
A retract of a contractible subshift is contractible.
\end{lemma}

\begin{proof}
Suppose $Y$ is contractible, and consider directly a subshift $X \subset Y$, and a retraction $r : Y \to X$ with $r|X = \id_X$. Let $h : I \times Y \times Y \to Y$ be the contraction homotopy. Define
\[ h'(t, x, x') = r(h(t, x, x')). \]
for $(t, x, x') \in I \times X \times X$. Observe that the image of $h(t, x, x')$ is in $Y$, so the composition makes sense, and we can take the codomain of $h'$ to be $X$. Now $h'(\bar 0, x, x') = r(x) = x$ since $h$ is a homotopy and $r$ a retraction, and similarly $h'(\bar 1, x, x') = x'$.
\end{proof}

\begin{theorem}
\label{thm:RetractCharacterization}
Let $X \subset \Sigma^G$ be a subshift. Then $X$ is a retract of $\Sigma^G$ (through its natural embedding, equivalently abstractly) if and only if all of the following hold:
\begin{itemize}
\item $X$ has a fixed-point,
\item $X$ is of finite type, and
\item $X$ is contractible.
\end{itemize}
\end{theorem}

\begin{proof}
The previous three lemmas show that the three properties hold for retracts of $\Sigma^G$.

Now suppose the properties hold. There exists a fixed point in $X$, so for some symbol $0 \in \Sigma$ we have $\bar 0 \in X$. Let $h : I \times X \times X \to X$ be the contraction homotopy. Consider the natural relaxation $\tilde h : I \times \Sigma^G \times \Sigma^G \to \Sigma^G$.

Fix a window $W \Subset G$ for $X$. Let $M \subset G$ be finite, and define $t_M : \Sigma^G \to I$ by
\[ t_M(x)_a = 0 \iff x|aM \mbox{ is globally valid in $X$}. \]
Define $r_M : \Sigma^G \to \Sigma^G$ by
\[ r_M(x) = \tilde h(t_M(x), x, \bar 0). \]
For $x \in X$, we have $r_M(x) = \tilde h(\bar 0, x, \bar 0) = x$, so it suffices to show that for a suitable choice of $M$, the codomain of $r_M$ can be restricted to $X$.

Suppose that, for some $M \subset G$, the codomain of $r_M$ cannot be restricted to $X$, \textit{i.e.} there are some configurations in the image of $r_M$ that do not belong to $X$. Then, using the fact $W$ is a window size for $X$ and possibly shifting, we find $x_M \in \Sigma^G$ such that
\[ \tilde h(t_M(x_M), x_M, \bar 0)|W \not\sqsubset X. \]

Note that since $0^W \sqsubset X$, we must have $\tilde h(t_M(x_M), x_M, \bar a)_b \neq 0$ for some $b \in W$. Thus by the construction of $\tilde h$, $0 \sqsubset t_M(x_M)|WN$ where $N$ is the neighborhood of $\tilde h$. By the definition of $t_M$, $x_M|bcM$ is globally valid in $X$ for some $b \in W, c \in N$. Note that $W$ and $N$ are fixed sets, while $M$ can be picked arbitrarily large, so for large $M$, $x_M$ will have arbitrarily large globally valid central patterns.

From this we conclude that 
we can have $\tilde h(t, x_M, \bar 0)|W \not\sqsubset X$ with $d(x_M, X)$ arbitrarily small. 
This means $\tilde h$ is not uniformly continuous in $I_2 \times X \times X$.  This contradicts the compactness of $I_2 \times X \times X$, continuity of $\tilde h$, and the fact $\tilde h$ is an extension of $h$.
\end{proof}


\subsection{Adaptation to general contractible subshifts}

It seems difficult to give a full characterization of retracts of subshifts. However, given that the characterization of retracts of full shifts is in terms of contractible subshifts, it is natural to characterize the retracts of contractible subshifts. As we see, the proof of this is essentially the same.




\begin{definition}[Subshift-relative contractibility]
Let $X \subset Y \subset A^G$ be subshifts. We say that $X$ is \emph{relatively contractible in $Y$} if there exists a block-map $h : I \times X \times X \rightarrow Y$ such that for all $x,y \in X$, $h\left(\bar 0,x,y\right) = x$, $h\left(\bar 1,x,y\right) = y$. 
\end{definition}

We call $h$ a \emph{relative contraction homotopy}. The definition says that, given configurations in $X$, $h$ only introduces patterns that are forbidden in $X$, but not ones forbidden in $Y$, and at endpoints it behaves like a standard contraction homotopy.

\begin{theorem}
\label{thm:GeneralRetractCharacterization}
Given two subshifts $X \subset Y$, if $Y$ is contractible, then $X$ is a retract of $Y$ if the following hold:
\begin{itemize}
\item $X$ is a relative SFT of $Y$,
\item there exists a map from $Y$ to $X$, and
\item $X$ is contractible. 
\end{itemize}
Conversely, the first two items hold whenever $X$ is retract of $Y$, and the third also holds if $X$ is additionally relatively contractible in $Y$.
\end{theorem}

If $Y$ is itself contractible, then any subshift of it is relatively contractible, as the restriction of the contraction homotopy has the required properties. Thus in this case the theorem provides a full characterization of retracts, generalizing Theorem~\ref{thm:RetractCharacterization} (noting that the existence of a morphism from a full shift to $X$ is equivalent to the existence of a fixed point).

\begin{corollary}
If $Y$ is contractible, then $X$ is a retract of $Y$ if and only if $X$ is a relative SFT of $Y$, there exists a map from $Y$ to $X$, and $X$ is itself contractible.
\end{corollary}


\begin{proof}[Proof of Theorem~\ref{thm:GeneralRetractCharacterization}]
Suppose first that the three properties hold. In other words,  let $X$ be a relative SFT of $Y$, let $f : Y \to X$ be a morphism, and let $h : I \times X \times X \to X$ be the contraction homotopy. We then construct a retract as in Theorem~\ref{thm:RetractCharacterization}: we take $\tilde h$ a natural relaxation of $h$, and the retraction is given by $r(x) = \tilde h(t(x), x, f(x))$ where $t(x)$ has large areas of $1$s where there are relative forbidden patterns of $X$ nearby, and $0$s elsewhere. We leave the details to the interested reader.

Conversely, if $r : Y \to X$ is a retraction, then
\begin{itemize}
\item $X$ is a relative SFT of $Y$ because it is the set of fixed points of $r$, and
\item $r : Y \to X$ is a morphism from $Y$ to $X$. 
\end{itemize}
For the third item, if $h' : I \times X \times X \to Y$ is a relative contraction homotopy for $X$ relative to $Y$, then $h = r \circ h'$ is a contraction homotopy for $X$.
\end{proof}

\section{Contractibility and homotopy equivalence}

In this section, we attempt to work out some basic ``algebraic topology'' for subshifts, based on the idea of a homotopy. We obtain the useful fact that in the case of contractible SFTs, homotopies can be composed. We also characterize contractibility of one-dimensional SFTs, and homotopy equivalence for transitive SFTs.


\begin{definition}
We say two subshifts $X, Y$ are \emph{homotopy equivalent}, and we write $X \cong Y$, if there are morphisms $f : X \to Y$, $g : Y \to X$ such that $g \circ f \cong \id_X$ and $f \circ g \cong \id_Y$.
\end{definition}

It is useful to allow auxiliary symbols in the parameter space of a homotopy:

\begin{definition}
\label{def:kSymbol}
A \emph{$k$-symbol homotopy} is $h : I_k \times X \to Y$ such that $h(\bar 0, x) = f(x), h(\overline{k-1}, x) = g(x)$ for all $x \in X$.
\end{definition}

It is very easy to see that this does not change the meaning of homotopy. We give a somewhat roundabout complicated proof, as it contains some ideas that are useful later.

For the proof, and also later in the text, we think of $I = I_2$ as a \emph{two-pointed subshift}, meaning it has a \emph{left endpoint} $\bar 0$ and a \emph{right endpoint} $\bar 1$, and these are taken to be part of the structure. We say a morphism between two two-pointed subshifts \emph{respects endpoints}, or is \emph{mod endpoints} if it maps the left (resp.\ right) endpoint to the left (resp.\ right) endpoint.

In general, for any two-pointed subshift $J$ we have a notion of \emph{J}-homotopy between $f$ and $g$, meaning 
$h : J \times X \to Y$ such that $h(y, x) = f(x), h(z, x) = g(x)$ for all $x \in X$, where $y, z$ are the two endpoints of $J$.


Let's say two two-pointed subshifts are \emph{homotopy equivalent mod endpoints} if we can pick the $f, g$ in the definition of homotopy equivalence so that they respect endpoints, and the homotopies between $g \circ f \cong \id_X$ and $f \circ g \cong \id_Y$ are mod endpoints meaning $h(t, x) = x$ for both endpoints $x$ and all $t$.

\begin{lemma}
\label{lem:mn}
$I_m$ and $I_n$ are homotopy equivalent mod endpoints for any $m,n \geq 2$.
\end{lemma}

\begin{proof}
Define $f : I_m \to I_n$ by
\[ f(x)_a = \left\{\begin{array}{ll}
0 & \mbox{if } x_a \neq m-1,  \\
n-1 & \mbox{otherwise}.
\end{array}\right.\]
and $g : I_m \to I_n$ by
\[g(x)_a = \left\{\begin{array}{ll}
0 & \mbox{if } x_a \neq n-1, \\
m-1 & \mbox{otherwise}.
\end{array}\right.. \]
Clearly these maps respect endpoints. The map $g \circ f$ preserves $(m-1)$-symbols, and maps everything else to $0$. For the homotopy $h : g \circ f \to \id_{I_m}$ we can use the naive homotopy:
\[ h(t, x)_a = \left\{\begin{array}{ll}
g(f(x))_a & \mbox{if } t_a = 0, \\
x_a & \mbox{if } t_a = 1.
\end{array}\right. \]
It is easy to check that $h(t,x) = x$ for all $t \in I_2$ when $x$ is an endpoint. The case of $f \circ g$ is symmetric.
\end{proof}

\begin{lemma}
For any $m,n \geq 2$, two morphisms are $m$-symbol homotopic if and only if they are $n$-symbol homotopic.
\end{lemma}

\begin{proof}
Suppose we have an $m$-symbol homotopy $h : I_m \times X \to Y$ between morphisms $f, g$. Let $\pi : I_m \to I_n$ be any morphism respecting endpoints. Precompose $h$ with $\pi \times \id_X$ to get $\hat{h} : I_n \times X \to Y$, i.e.\ $\hat{h}(t, x) = h(\pi(t), x)$. Then $\hat{h}(\bar 0, x) = h(\bar 0, x) = f(x)$ and $\hat{h}(\overline{m-1}, x) = h(\overline{n-1}, x) = g(x)$.
\end{proof}

Two two-pointed subshifts can be joined in an obvious way by identifying the right endpoint of the first with the left endpoint of the second, analogously to the usual definition of path composition when defining the fundamental group in algebraic topology. In our symbolic category, this turns out a little more subtle, since already the composition of two copies of $I_2$, is not isomorphic to $I_2$. Indeed, it is not even a subshift of finite type. We will see that, nevertheless, in the category of subshifts of finite type, the resulting generalized notion of homotopy does not change.

Denote by $I_{a,b}$ the subshift isomorphic to $I_2$ on symbols $\{a, b\}$, with left and right endpoints $\bar a, \bar b$.

\begin{definition}
\label{def:Pk}
Let $k \geq 1$. The \emph{$k$-path shift} is $P_k = I_{0,1} \cup I_{1,2} \cup \cdots \cup I_{k-1,k}$. A \emph{$k$-step homotopy} is $h : P_k \times X \to Y$ such that $h(\bar 0, x) = f(x), h(\bar k, x) = g(x)$ for all $x \in X$. A \emph{multistep homotopy} is a $k$-step homotopy for some $k$.
\end{definition}

Note that $P_k$ is naturally a two-pointed subshift with endpoints $\bar 0, \bar k$. Thus a $k$-step homotopy is simply a $P_k$-homotopy in the terminology introduced below Definition~\ref{def:kSymbol}. We have $P_1 = I_2$, and $P_k$ is of finite type if and only if $k = 1$. 


\begin{lemma}
\label{lem:StarComposable}
For subshifts $X, Y$, multistep homotopy is an equivalence relation on morphisms $f : X \to Y$.
\end{lemma}

\begin{proof}
One can invert and compose homotopies in an obvious way.
\end{proof}

\begin{lemma}
Every SFT approximation of $P_m$ is homotopy equivalent to $I_2$ mod endpoints.
\end{lemma}

\begin{proof}
An SFT approximation of $P_m$ simply checks that the symbols at $a, b$ differ by at most $1$ if $a$ and $b$ belong to a translate of a particular finite subset of the group. The idea is then that for the map from $I_2$ to $P_m$, if we are deep inside an area of $0$s, then we write $0$, and otherwise we write a naive transient towards the maximal symbol, depending on the distance from an area of $0$s.

In formulas, define $\delta : G \times I_2 \to \N \cup \{\infty\}$ by $\delta(a, x) = \inf \{d_G(a, b) \;|\; b \in G, x_b = 0\}$. Then for some large $r$ define $f : I_2 \to P_m$ by
\[ f(x)_a = \min(\lfloor \delta(a, x)/r \rfloor, m). \]
We claim that $f$ is well-defined, i.e.\ if $r$ is sufficiently large, then the image is contained in any given SFT approximation of $P_m$. For this, simply note that $\delta(a, x)/r$ is $\frac{1}{r}$-Lipschitz in $a$ for a fixed $x$, and combine this with the first sentence of the proof.

For $g : P_m \to I_2$ we use simply the map
\[ g(x)_a = \min(x_a, 1). \]

Note that $g \circ f$ and $f \circ g$ indeed respect endpoints. The map $g \circ f$ is homotopy-equivalent to $\id_{I_2}$ by the naive homotopy, which is mod endpoints. Now consider the map $f \circ g$. What this map does is it forgets the positions of all non-zero symbols, and replaces them by a (discrete) speed-$\frac1r$ transient depending on the distance from the nearest $0$-symbol, terminating at the $m$-symbol once far enough.

We can construct a homotopy $h : \id_{P_m} \to f \circ g$ analogously to the above definition of $f$: define $L : [0,1] \times \{0, \ldots, m+1\}^2 \to \{0, \ldots, m+1\}$ by $L(t, i, j) = \lfloor i + (j - i) r \rfloor$ (one can think of this as a discrete analog of a linear interpolation map -- or as a discrete homotopy). Then define
\[ h(t, x)_a = L(\min(\delta(a, t)/r, 1), x_a, f(g(x))_a). \]
The calculation $h(\bar 0, x)_a = L(0, x_a, f(g(x))_a) = x_a$ for all $x \in P_m, a \in G$ and $h(\bar 1, x)_a = L(1, x_a, f \circ g) = f(g(x))_a$ shows that this is a homotopy, and it respects endpoints because $f,g$ respect endpoints and $L(t, i, i) = i$ for all $t, i$.
\end{proof}

The reader may notice that in the previous proof, 
no properties of the specific morphism $f \circ g$ were used to define the homotopy (other than checking the preservation of endpoints). Indeed, any two endomorphisms are homotopic on this subshift, by the same formula. This reflects the fact that SFT approximations of $P_m$ are contractible subshifts, as we will see.

\begin{lemma}
\label{lem:Transitive}
If $Y$ is a subshift of finite type, then $k$-step homotopy is equivalent to $1$-step homotopy for morphisms $f : X \to Y$.
\end{lemma}

\begin{proof}
Consider a $k$-symbol homotopy $h : I_k \times X \to Y$ from $f$ to $g$ (recall that $k$-symbol homotopy is equivalent to $2$-symbol homotopy). Let $\pi : P_k \to I_k$ be the natural inclusion. Then $h \circ (\pi \times \id_X) : P_k \times X \to Y$ is a homotopy from $f$ to $g$, thus is a $k$-step homotopy.

Consider then a $k$-step homotopy $h : P_k \times X \to Y$ from $f$ to $g$. Since $Y$ is a subshift of finite type, there is an SFT approximation $P$ of $P_k$ such that a domain extension $\tilde h : P \times X \to Y$ of $h$ exists. Let $\pi$ be a morphism from $I_2$ to $P$ which respects endpoints (endpoints of $P$ being those of $P_k$) from the previous lemma. Then $\tilde h \circ (\pi \times \id_X) : I_2 \times X \to Y$ is a homotopy from $f$ to $g$.
\end{proof}

\begin{lemma}
Let $X$ be a subshift containing a fixed point $\bar a$, let $\pi_1, \pi_2 : X \times X \to X$ be the projections. Then
\[ \id_X \cong_k (x \mapsto \bar a) \implies \pi_1 \cong_{2k} \pi_2 \]
and
\[ \pi_1 \cong_k \pi_2 \implies \id_X \cong_k (x \mapsto \bar a). \]
\end{lemma}

\begin{proof}
We simply mimic the proof of Lemma~\ref{lem:TwoHomotopyDefinitions}, though in the first case we have to replace the explicit formula by an abstract composition of two homotopies. More precisely, suppose $h : \id_X \cong_k (x \mapsto \bar a)$. Then 
\[ h'(t, x, y) = h(t, x) \]
is a homotopy from $\pi_1$ to $((x, y) \mapsto \bar a)$, and
\[ h'(t, x, y) = h(\bar 1 - t, y) \]
is a homotopy from $((x, y) \mapsto \bar a)$ to $\pi_2$. The proof of Lemma~\ref{lem:StarComposable} shows that $\pi_1 \cong_{2k} \pi_2$.

In the second case we can directly copy the formula from Lemma~\ref{lem:TwoHomotopyDefinitions}: if $h'$ is a homotopy between the two projections, then
\[ h(t, x) = h'(t, x, \bar 0) \]
gives a homotopy between the identity map on $X$, and the map from $X$ to $\{\bar 0\}$.
\end{proof}

\begin{corollary}
\label{cor:FixedPointContractible}
If $X$ is a subshift of finite type with a fixed point, then it is fixed-point contractible if and only if it is contractible. In particular, in this case fixed-point contractibility does not depend on the fixed-point.
\end{corollary}

This indeed needs the finite type assumption, as $P_2 = I_{0,1} \cup I_{1, 2}$ is fixed-point contractible to $\bar 1$, but not $\bar 0$, and is not contractible. Of course, by the above lemma, we do have a $2$-step homotopy $\pi_1 \cong_2 \pi_2$ between the two projections $\pi_1, \pi_2 : P_2 \times P_2 \to P_2$. 

\begin{lemma}
Homotopy equivalence is an equivalence relation for subshifts of finite type.
\end{lemma}

\begin{proof}
Suppose $X, Y, Z$ are subshifts of finite type, and suppose $X \cong Y \cong Z$. There exist $f_1 : X \to Y, f_2 : Y \to Z, g_1 : Y \to X, g_2 : Z \to Y$ such that $g_1 \circ f_1 \cong \id_X, f_1 \circ g_1 \cong \id_Y, g_2 \circ f_2 \cong \id_Y, f_2 \circ g_2 \cong \id_Z$.

Define $f = f_2 \circ f_1 : X \to Z, g : g_1 \circ g_2 : Z \to X$. We claim that these give a homotopy equivalence. Consider $g \circ f : X \to X$. This is just the map $g_1 \circ g_2 \circ f_2 \circ f_1$. We first form a homotopy from $g_1 \circ g_2 \circ f_2 \circ f_1$ to $g_1 \circ f_1$. For this, let $h_2 : g_2 \circ f_2 \cong \id_Y$. We define the block map
\[ h(t, x) = g_1(h_2(t, f_1(x))) \]
and observe that
\[ h(\bar 0, x) = g_1(h_2(\bar 0, f_1(x))) = g_1(g_2 \circ f_2(f_1(x)))) \]
and
\[ h(\bar 1, x) = g_1(h_2(\bar 1, f_1(x))) = g_1(\id_Y(f_1(x)))) = g_1(f_1(x)) \]
as desired. From $g_1 \circ f_1 \cong \id_X$ it follows from Lemma~\ref{lem:StarComposable} that $g \circ f$ and $\id_X$ are $2$-step homotopic, and since $X$ is of finite type, they are homotopic. Similarly, one shows that $f \circ g \cong \id_Y$.
\end{proof}

\begin{theorem}
\label{thm:HEToPoint}
Let $X$ be a subshift of finite type. Then $X$ is homotopy equivalent to a one-point subshift if and only if it is contractible and has a fixed-point.
\end{theorem}

\begin{proof}
If $X$ does not have a fixed-point, then it of course cannot be homotopy equivalent to a one-point subshift (as there is no map from the one-point subshift into $X$). So we simply need to show that for subshifts of finite type with a fixed-point, homotopy equivalence with a one-point subshift is equivalent to contractibility.

Suppose $\bar 0 \in X$. Suppose first that $X$ is contractible. %
By the previous corollary, it is fixed-point contractible to $\bar 0$, meaning $h : \id_X \cong \bar 0$. Let $f : X \to \{\bar 0\}$ be the constant map, and let $g : \{\bar 0\} \to X$ be the inclusion morphism. Then $g \circ f$ is the endomorphism of $X$ that maps all points to $\bar 0$. By assumption, we have $h^R : g \circ f \cong \id_X$, where $h^R$ is the inverse homotopy obtained by changing the rules of $0$ and $1$. On the other hand, $f \circ g = \id_{\{\bar 0\}}$ so certainly these maps are homotopic. Thus, $f$ and $g$ form a homotopy equivalence.

Next suppose that $X$ is homotopy equivalent to a one-point subshift. We may assume this subshift is $\bar a$. It follows that there exists maps $f : X \to \{\bar a\}$ and $g : \{\bar a\} \to X$ such that $g \circ f \cong \id_X$ and $f \circ g \cong \id_{\{\bar a\}}$. Now $g \circ f$ is an endomorphism of $X$ that maps all points to some $\bar b = g(\bar a)$. So already $g \circ f \cong \id_X$ implies that $X$ is fixed-point contractible, thus contractible.
\end{proof}

\begin{lemma}
If $X$ is contractible and $X$ is homotopy equivalent to $Y$, then $Y$ is $3$-step contractible. In particular, if $Y$ is of finite type, it is contractible.
\end{lemma}

\begin{proof}
This is clear from the previous lemma if $X$ is a subshift of finite type with a fixed point. We give a direct proof for the general case. The assumptions are that there is a homotopy $h : \pi_1 \cong \pi_2$ where $\pi_1, \pi_2 : X \times X \to X$ are the two projections, and that there exist $f : X \to Y$ and $g : Y \to X$ such that $g \circ f \cong \id_X$ and $h' : f \circ g \cong \id_Y$.

Concretely $h' : f \circ g \cong \id_Y$ means
\[ h' : I_2 \times Y \to Y \]
such that $h'(\bar 0, y) = f \circ g(y)$ and $h'(\bar 1, y) = y$.

We define
\[ h''(t, y, y') = f(h(t, g(y), g(y'))) : I_2 \times Y \times Y \to Y. \]
Clearly this is indeed a morphism (and the codomain choice is valid). We have
\[ h''(\bar 0, y, y') = f(h(\bar 0, g(y), g(y'))) = f(g(y)) \]
and
\[ h''(\bar 1, y, y') = f(h(\bar 1, g(y), g(y'))) = f(g(y')) \]
so this is a homotopy from $f \circ g \circ \pi_1$ to $f \circ g \circ \pi_2$.

Next we observe that $\pi_1 \cong f \circ g \circ \pi_1$: simply use $h'''(t, y, y') = h'(\bar 1 - t, y)$ (where subtraction is performed cellwise, $(z - z')_a = z_a - z'_a$). Similarly, $f \circ g \circ \pi_2 \cong \pi_2$. We conclude $\pi_1 \cong_2 \pi_2$, i.e.\ $3$-step contractibility. Finally, when $Y$ is SFT, we obtain contractibility by Lemma~\ref{lem:Transitive}.
\end{proof}

We now conclude that $P_m$ has contractible SFT approximations, as was claimed previously.

\begin{proposition}
\label{prop:PmSFTApproContractible}
Every SFT approximation of $P_m$ is contractible.
\end{proposition}

\begin{proof}
We showed above that every SFT approximation of $P_m$ is homotopy equivalent to $I_2$ (even mod endpoints). Since $I_2$ is contractible and SFT approximations are SFT, we conclude from the previous lemma that indeed SFT approximations of $P_m$ are contractible.
\end{proof}

\begin{lemma}
If $Y$ is contractible, then any two morphisms $f,g : X \to Y$ are homotopic.
\end{lemma}

\begin{proof}
Let $h : I_2 \times Y \times Y \to Y$ be a homotopy from $\pi_1$ to $\pi_2$. Map $h'(t, y) = h(t, f(y), g(y))$.
\end{proof}

The following very simple ``time dilation lemma'' shows that in the time parameter $t$, we may always assume all $0$-areas are contained in large balls. It would be possible to simultaneously say something nice about $1$-areas, at least on abelian groups, but we do not know in which generality one can make both types of areas consist of literal balls. 

\begin{lemma}
\label{lem:TimeDilation}
Let $n \geq 1$ be arbitrary. Let $J \subset I_2$ be the subshift (depending on $n$) where for each $x \in J$, we require that $x_a = 0 \implies x|bB_n \equiv 0$ for some $b \in aB_n$. Then the notion of $J$-homotopy (with the same endpoints $\bar 0, \bar 1$) is equivalent to the notion of $I$-homotopy.
\end{lemma}

\begin{proof}
An $I_2$-homotopy immediately gives a $J$-homotopy by restriction. From a $J$-homotopy, we obtain an $I_2$-homotopy by constructing an endpoint-respecting morphism from $I_2$ to $J$.

Define $f : I_2 \to I_2$ by $f(t)_a = \min(\{t_b \;|\; b \in a B_n\})$. This clearly fixes the endpoints. For $t \in I$, whenever $f(t)_a = 0$, we have $0 \sqsubset t|aB_n$, say $t_{ab} = 0$ for $b \in B_n$. Then of course $f(t)|abB_n \equiv 0$, and we conclude that $f(t) \in J$. 
\end{proof}

For the group $\Z$, we can easily turn both $0$- and $1$-areas into long intervals.

\begin{lemma}
\label{lem:TimeDilationZ}
For $G = \Z$, let $J \subset I_2$ be the subshift where for each $x \in J$  $b \in \{0,1\}$ and $a \in \Z$, we require that $x_a = b \implies x|[i, i+n] \equiv b$ for some $i \in \Z$ such that $a \in [i, i+n]$. Then the notion of $J$-homotopy (with the same endpoints $\bar 0, \bar 1$) is equivalent to the notion of $I_2$-homotopy.
\end{lemma}

\begin{proof}
Again it suffices to construct an endpoint-respecting morphism from $I_2$ to $J$. For this, apply first $f(t)_a = \min(t|[a, a+10n])$ (an endomorphism of $I_2$) and then $g(t)_a = \max(t|[a, a+n])$ (another endomorphism of $I_2$), and observe that the image is always in $J$, and endpoints are respected.
\end{proof}

\begin{theorem}
\label{thm:MixingContractible}
In one dimension, an SFT is contractible if and only if it is mixing.
\end{theorem}

\begin{proof}
If a subshift is contractible, it is strongly irreducible, thus mixing.

To show that a mixing $\Z$-SFT $X$ is contractible, observe first that contractibility is conjugacy-invariant, so we may assume $X$ has window $\{0,1\}$. Recall \cite{LiMa21} that mixing is equivalent to the existence of $n$ such that for any two symbols $a, b$, there exists $v = v(a, b)$ such that $|v| = n$ and $avb \in L(X)$. Let us index the word $v$ as $v_0 v_1 \ldots v_{n-1}$.

By the previous lemma, it suffices to construct a $J$-homotopy between the two projections, where in $J$ every $0$- and $1$-symbol belongs to an interval of length at least $n+1$. Define $h : J \times X \times X \to X$ by
\[ h(t, x, y)_i = \left\{\begin{array}{ll}
x_i & \mbox{if } t|[i, i+n] \equiv 0 \\
v(x_{i+j-n-1}, y_{i+j})_{n-j} & \mbox{if } \exists j \in [1, n]: t|[i, i+j-1] \equiv 0, t_{i+j} = 1 \\
y_i & \mbox{if } t|[i, i+n] \equiv 1 \\
v(y_{i+j-n-1}, x_{i+j})_{n-j} & \mbox{if } \exists j \in [1, n]: t|[i, i+j-1] \equiv 1, t_{i+j} = 0 \\
\end{array}\right. \]

Since we restrict to configurations in $J$, it is straightforward to verify that this is a homotopy between the two projections, as desired, namely between coordinates $i < j$ where we copy from different points, there is at some point a change in $t$, and at the change we introduce a valid transient using $v$.
\end{proof}

\begin{theorem}
\label{thm:HomotopyEquivalent}
In one dimension, two transitive SFTs are homotopy equivalent if and only if they map into each other.
\end{theorem}

\begin{proof}
If there is no map from $X$ to $Y$ (or vice versa), then $X$ and $Y$ clearly cannot be homotopy equivalent. We show that they are homotopy equivalent in all other cases.

For mixing SFTs, the claim is immediate from the previous theorem, namely if we have maps $f : X \to Y$ and $g : Y \to X$, by contractibility of $X$ and $Y$, and by the previous lemma, we have $f \circ g \cong \id_Y$ and $g \circ f \cong \id_X$.

To deal with general transitive SFTs, we recall basic structure theory of SFTs (see \cite{LiMa21}). Namely, we consider the \emph{primitive period} of a transitive SFT, which is the least $p \in \N_+$ such that we can find a clopen cross-section $C$, meaning that every $x \in C$ returns to $C$ in exactly $p$ steps, and the first-return map to $C$ is mixing. Each $x \in X$ then has a well-defined \emph{phase} $i \in \Z_p$.

The primitive period cannot increase under a morphism (by mixing of the first-return map), so two transitive SFTs that map into each other must have the same primitive period.

Consider now two transitive SFTs $X, Y$ which map into each other. These maps $f : X \to Y, g : Y \to X$ must act consistently on the phase, in the sense that if some $x \in X$ (resp.\ $y \in Y$) has phases $i$ and $f(x)$ (resp.\ $g(y)$) has phase $j$, then every point with phase $i'$ is mapped to a point with phase $i' + (j-i)$. Thus by composing $f$ and $g$ with shifts we may assume they map the chosen cross-sections of $X$ and $Y$ into each other.

The first-return map is well-known to be of finite type, so as a mixing SFT it is contractible by the above proof. It follows that the maps $f, g$ give a homotopy equivalence between the first-return maps of $X, Y$ to their chosen cross-sections, by the first part of the proof. The homotopy must automatically respect phases of points, so in fact $f, g$ are a homotopy equivalence between $X$ and $Y$.
\end{proof}

Two transitive $\Z$-SFTs $X$ and $Y$ map into each other if and only if $P(X) = P(Y)$ where $P(Z)$ is the set of (not necessarily least) periods of points in a subshift $Z$ \cite{LiMa21}. 

\section{Periodic asymptotic dimension, almost unions, and dense periodic points} 

In this section, we define the periodic asymptotic dimension (of a group) and the patching property (of a subshift). We show that if a subshift has the patching property and the acting group $G$ has finite periodic asymptotic dimension, then periodic points are dense. We show that contractibility and FEP imply the patching property, so we obtain dense periodic points for these subshifts. Note that in Section~\ref{sec:Patching} we also define a patching property for groups.


\subsection{Periodic asymptotic dimension}

Periodic asymptotic dimension is a variant Gromov's asymptotic dimension. To define it, we define the latter in a symbolic dynamical way so that it corresponds to the nonemptiness of certain subshifts. Our dimension is then obtained by replacing nonemptiness with existence of periodic points. Let $G$ be a group.

\begin{definition}
Let $X_{G, n, r, m}$ (or just $X_{n, r, m}$ when dealing with a single group $G$) be the subshift of the full shift $\{1,\ldots,n\}^G$, containing those configurations $x \in \{1,\ldots,n\}^G$ such that for each $1 \leq c \leq n$, if we define a graph with nodes $\{a \in G \;|\; x_a = c\}$ and edges $\{(a, b) \in G^2 \;|\; x_a = x_b = c \wedge |a^{-1}b| \leq r\}$, then the connected components of the graph are of cardinality at most $m$.
\end{definition}

Note that $X_{G, n, r, m}$ is a subshift of finite type for any $G, n, r, m$. In words, in $X_{G, n, r, m}$ we require that for each color $1 \leq c \leq n$, the graph with $c$-colored nodes and edges between $r$-separated nodes has components of cardinality at most $m$. If $x \in X_{n, r, m}$, we call the sets $x^{-1}(i) \subset G$ \emph{color classes}, and refer to the (bounded) components of the aforementioned graph \emph{(monochromatic) $r$-components}.

We usually think in terms of $r$-paths, meaning sequences $a_1, a_2, \ldots$ in $G$ such that $d(a_i, a_{i+1}) \leq r$ for all $i$, and the assumption that $x \in X_{n, r, m}$ for some $m$ is equivalent to showing that there is a uniform bound on the lengths of such injective monochromatic paths. Often we also say there is an \emph{$r$-jump} if $d(a, b) \leq r$, so an $r$-path is sequence of $r$-jumps.

Recall that one of the equivalent definitions of Gromov's asymptotic dimension is that it is the infimum of all $d$ such that $\forall r: \exists m: X_{d+1, r, m}$ is nonempty.

We define a periodic variant of this:

\begin{definition}
The \emph{periodic asymptotic dimension} $\pdim(G)$ of a group $G$ is the infimum of all $d$ such that 
\[ \forall r: \exists m: X_{d+1, r, m} \mbox{ has a periodic point}. \]
\end{definition}

(Recall that a periodic point is one having finite shift-orbit.)

We can of course generalize asymptotic dimension and periodic asymptotic dimension to a general countable group $G$ by instead of word norm picking sets $S_r \Subset G$, where $G = \bigcup_r S_r$ as an increasing union, and by $r$-paths $p : I \to G$ referring to paths with offsets $p(i)^{-1} p(i+1) \in S_r$. It is easy to see that the dimension obtained does not depend on the choice of the sequence $(S_r)_r$. 

Recall that the \emph{residual finiteness core} \cite{BrJa23} of a group $G$ is the intersection of all of its finite-index subgroups.

\begin{lemma}
\label{lem:RFC}
If an infinite f.g.\ group $G$ has finite periodic asymptotic dimension, then its residual finiteness core is locally finite.
In particular, if $G$ is torsion-free, then it must be residually finite.
\end{lemma}

\begin{proof}
Suppose $G$ has periodic asymptotic dimension $d$. Suppose that the residual finiteness core contains an infinite finitely-generated subgroup $H$, generated by say $T \Subset H$. Pick $r$ such that $T \subset B_r$. We claim that there does not exist $m$ such that $X_{d+1, r, m}$ contains a periodic point.

Suppose on the contrary that there is a periodic point $x \in X_{d+1, r, m}$ of $G$ for some $m$. Since the stabilizer $K$ of this configuration contains $H$, we have $x|H = i^H$ where $i = x_{1_G}$. Since $H$ is generated by $T$, this is contained in a single monochromatic component, contradicting $x \in X_{d+1, r, m}$.
\end{proof}

We show in Section~\ref{sec:Dimension} for example that all virtually polycyclic groups have finite periodic asymptotic dimension, and so do metabelian Baumslag-Solitar groups and the lamplighter group $\Z_2 \wr \Z$. Furthermore, the dimension of $\Z^d$ is the expected $d$ for this notion.

There are many equivalent definitions of the asymptotic dimension. We will need the analog of one of them. 
Namely, we define a subshift $Y_{d+1, r, m}$ which, instead of recording actual color classes, only records their shapes. This can be coded into a finite alphabet by explicitly recording the relative shape of an $r$-component containing a cell into its color. The following definition makes this explicit.

\begin{definition}
Let $B = S^{\leq rm}$, and take as alphabet $\Sigma$ of $Y_{n, r, m}$ the power set $P(B)$. Now, any $y \in \Sigma^G$ %
can be interpreted as a graph with vertices $G$ and with a directed edge from $g$ to $gh$ whenever $h \in y_g$. Now we define $Y_{n, r, m}$ by requiring that
\begin{itemize}
    \item this graph is a union of cliques;
    \item the components (i.e.\ the cliques) are $r$-connected as subsets of $G$;
    \item the components are of cardinality at most $m$;
    \item any $r$-ball in $G$ intersects at most $n$ distinct components.
\end{itemize}
\end{definition}

Note that the first item in practice is the SFT condition $gx_g = h \implies h^{-1} \in x_{gh}$, the second and third are restrictions on the actual alphabet allowed in $Y_{n, r, m}$, and the fourth one concretely says that for all $g \in G$, we can partition $gS$ into at most $n$ pieces $A_1, A_2, ..., A_{d+1}$ so that for all $i$ and all $g_1, g_2 \in A_i$ we have $g_2 \in g_1y_{g_1}$ (i.e. each $A_i$ is a subset of the same component).


\begin{lemma}
\label{lem:XYFactors}
For any $d$, there is a morphism $\phi : X_{d+1, 2r, m} \to Y_{d+1, r, m}$.
Conversely, if $r'$ is sufficiently large compared to $r$, and $m'$ is sufficiently large compared to $m$, then there is a morphism $\phi : Y_{d+1, r', m} \to X_{d+1, r, m'}$.
\end{lemma}

\begin{proof}
For the first claim, we can simply write down the arrows to other cells in the color component of a cell (w.r.t.\ $r$-paths). Since the monochromatic $2r$-components in configurations of $X_{d+1, 2r, m}$ have cardinality at most $m$, certainly the same cardinality bound holds when we consider the monochromatic $r$-components of such a configuration. A ball of radius $r$ cannot touch more than two distinct monochromatic $2r$-components of the same color in a point of $X_{d,2r,m}$, by the definition of a monochromatic component.

Assume now that $r' > r(d+2)$. If $y \in Y_{d+1, r', m}$, we construct a coloring of $G$ with $\{1, \ldots, d+1\}$ by a shift-commuting procedure, and show that the $r$-components are bounded (by some function of $r', m$). The bound can be taken as $m'$.

For this, for cell $a$ we pick its color $\phi(y)_a$ to be the maximal $1 \leq k \leq d+1$ such that the $rk$-ball around it touches $k$ or more components of $y \in Y_{d+1, r', m}$ (so it actually has to touch exactly $k$, or we would not reach the maximum, unless $k = d+1$).

We show that in the resulting configuration $x \in \{1, \ldots, d+1\}^G$, monochromatic $r$-components are bounded. Namely, suppose that a path with $r$-jumps in stays in a single color class of $\phi(y)$. Then we claim it has to actually stay near the same set of components of $y$. The reason is simply that if the color stays $k$ and some component is exchanged for another, the one dropped is now at distance at most $r(k+1)$, so in fact we should in fact have used the color $k+1$ (or larger) in the previous step, assuming $k+1 \leq d+1$.

But in fact, if $k = d+1$, then the $r(k + 1)$-ball still only touches at most $k$ color components in $y$, since $r' > r(d + 2)$. Thus, the logic also works for the maximal color.
\end{proof}

We use the same idea in Theorem~\ref{thm:ZdDimension}.

Note that by the lemma, we can also use $Y_{d+1,r,m}$ to characterize periodic asymptotic dimension, since morphisms preserve periods.

\begin{definition}
\label{def:Net}
Say that a subset $A$ of a metric space $B$ is an \emph{$R$-packing} if $\min \{d(a, b) \;|\; a, b \in A, a \neq b\} > R$, an \emph{$R$-covering} if $\max_{b \in B} d(b, A) \leq R$, and an \emph{$R$-net} if it is both.
\end{definition}

More generally, we say $A \subset B$ is an \emph{$R$-covering} of $C \subset B$ if $\max_{c \in C} d(c, A) \leq R$.


\begin{lemma}
\label{lem:ContractibleDimensionSFT}
On the group $\Z^d$, there exists $D$ such that for all $r$, the subshift $Y_{D, r, m}$ contains a nonempty contractible SFT for large enough $m$.
\end{lemma}

We can surely pick $D = d+1$, but the proof is geometrically simpler with a larger $D$.

\begin{proof}
Pick $D = 24^d$. We use the $\ell_\infty$ metric on $\Z^d$. To an $r$-net of vectors with integer coordinates, we associate a partition of the plane, and vice versa. There are many ways to do this, the most obvious possibility being Voronoi cells. However, the non-uniqueness of centers poses some technical difficulties, so we use a different method.

For each point $\vec v \in \Z^d$ consider the cone of all points whose coordinates are at least as large as those in the corresponding coordinate of the \emph{corner point} $\vec v$. For each $\vec u \in \Z^d$, associate $\vec u$ to the component of $\vec v$ if $\vec v$ is the closest corner point whose cone contains $\vec u$. We use the lexicographic order to break distance ties. It is then easy to see that components are $1$-connected, since if $\vec u$ is associated to the cone of corner $\vec v$, we stay in the same component when moving toward $\vec v$ along standard generators.

Clearly the sets of components are in shift-commuting continuous bijection with the corner points, so we see that this set of partitions assocated to $R$-nets is an SFT $Y$ (we simply check the $r$-net condition for the corner points).

It is easy to see that there is a bound on the radius of the components: If $\vec u \in \Z^d$, then there is a corner point $\vec v$ in the hypercube $\prod_i [\vec u_i, \vec u_i + 2r]$ by $r$-density, thus at distance at most $2r$. Thus, the associated closest corner point of $\vec u$ is also at distance at most $2r$. It follows that for each corner point $\vec v$, the set of associated vectors is a connected subset of $\prod_i [\vec v_i, \vec v_i + 2r]$, in particular its diameter is at most $2r$ as a subset of $\Z^d$.

On the other hand, if $\vec v$ is a corner point, then certainly all of the points $\vec u \in \vec v + [0, r/2]^d$ are associated to it, since we can have at most one corner point $\vec v$ such that $d(\vec u, \vec v) \leq r/2$, by the $r$-packing condition. Since $\lfloor r/3 \rfloor \leq r/2$, at least $r/3$ vectors are associated to $\vec v$.

Thus, the number of points associated to a single corner point ís a connected subset of $\Z^d$ whose size is at least $(r/3)^d$ and has diameter at most $2r$. Suppose now that $\vec u + B_r$ intersects more than $D$ distinct components. Then in the ball $\vec u + B_{3r}$ we have at least 
\[ D(r/3)^d = 24^d r^d / 3^d \geq (8r)^d \]
distinct elements. But $|\vec u + B_{3r}| = (6r + 1)^d < (8r)^d$, a contradiction.

We have shown that the SFT $Y$ (or rather, the encoding of it by writing in each node the edges to other elements of the component) is contained in $Y_{D, r, m}$ for any $m \geq (2r+1)^d$. Namely, the graph representation is by definition a union of cliques; its components are $1$-connected (thus $r$-connected); its components are of cardinality at most $(2r+1)^d \leq m$; any $r$-ball in $G$ intersects at most $D$ distinct components.



We claim that the SFT $Y$ is contractible. 
For this, it suffices to pick where the contraction homotopy writes corner points, given a time parameter $t$, and the sets of corner points of $x, y \in Y$. Then we simply write back the corresponding partition to get a point of $Y$.

In the time parameter $t$, when $t$ has a large area of $0$s, we of course want to copy from $x$, Specifically, in such areas, we pick corners from the first input point. We do the same for areas of $1$s. 

After this, we add all corner points from $y$ which are not at distance less than $r$ from points already included. After this step, we have a set of corner points which are an $r$-packing and a $2r$-covering. For the latter claim, any $\vec u$ will be at distance $r$ from a corner point $\vec v$ from $y$, and because this point was not included, there is an already-included point at distance $r$ from $\vec v$.

On an abelian group,\footnote{The argument specifically needs a bi-invariant metric. Since only virtually abelian groups have such metrics, the lemma does not immediately generalize to an essentially larger class of groups.} there is a general method for turning an $R$-covering $r$-packing $V \subset \Z^d$ into an $r$-net, in a shift-invariant continuous way.

This is just a greedy version of the usual existential proof that such nets exist. Namely, start with $V_0 = V$. List all vectors of length at most $R$, and when considering the vector $\vec v_i \in B_R$, set
\[ V_i = V_{i-1} \cup \{\vec u + \vec v_i \;|\; \vec u \in V, d(\vec u + \vec v_i, V_{i-1}) > r \}. \]

By induction, $V_i$ is an $r$-packing for all $i$: the added vectors are never in conflict with previously added ones (by definition of $V_i$), so the only possible conflict is that $d(\vec u + \vec v_i, \vec v + \vec v_i) \leq r$ for some $\vec u, \vec v \in V$. But $d(\vec u + \vec v_i, \vec v + \vec v_i) = d(\vec u, \vec v) > r$ by the assumption on $V$.

On the other hand $V_\ell$ for $\ell = |B_R|$ is an $r$-covering: if $\vec v \in \Z^d$, then there exists $\vec u \in V$ with $d(\vec v, \vec u) \leq R$. Thus, for some $i$ we have $\vec v_i = \vec v - \vec u \in B_R$. If $\vec u + \vec v_i = \vec v$ is not in $V_i$, then because it was not added to $V_i$ at step $i$, we must have $d(\vec u + \vec v_i, V_{i-1}) = d(\vec v, V_{i-1}) \leq r$, in particular $d(\vec v, V_\ell) \leq r$.
\end{proof}

\begin{remark}
It is possible to extend the definition of periodic asymptotic dimension so that the dimension is allowed to grow with scale, in the sense that if $d : \N \to \N$ is a function, we can say $G$ \emph{admits periodic asymptotic dimension} $d$ if
\[ \forall r: \exists m: X_{d(r)+1, r, m} \mbox{ has a finite subsystem}. \]
Our proofs should go through with $d$ growing sufficiently slowly rather than being constant.
\end{remark}

\subsection{Almost-unions of sets and patterns}

If $A \subset G$, write $A^{\circ r} = \{a \in G \;|\; aB_r \subset A\}$, where $B_r$ is the ball of radius $r$ (around the identity element, with respect to the fixed word metric).

\begin{definition}
\label{def:AlmostUnion}
Let $G$ be a finitely-generated group. The \emph{$r$-almost-union} of $A,B \subset G$ is defined by
\[ A \overset{r}{\cup} B = A^{\circ r} \cup B^{\circ r} \cup \left((A \cup B) \setminus (A B_r \cap B B_r)\right). \]
\end{definition}

\begin{figure}
\centering
\includegraphics[scale=0.61,trim=0 1.5cm 0 1cm,clip]{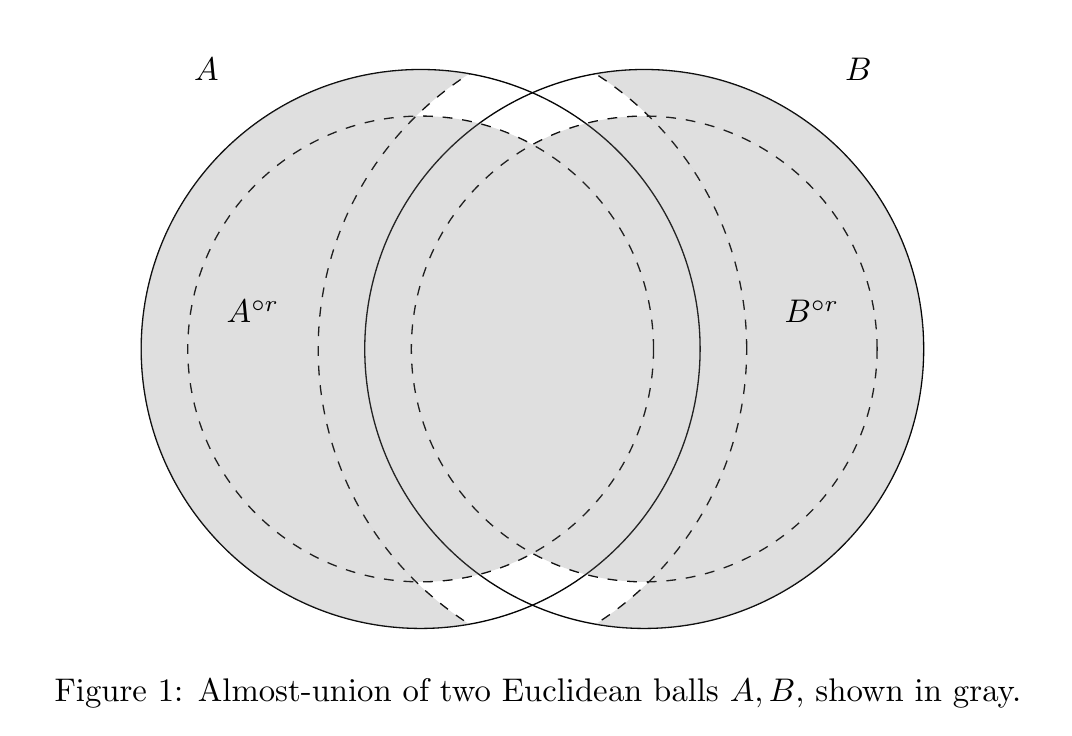}
\caption{Almost-union of two Euclidean balls $A, B$, shown in gray.}
\end{figure}

\begin{definition}
Let $X \subset \Sigma^G$ be a contractible subshift with gluing morphism $h$ with minimal radius $r$, and let $A, B \subset G$. Let $x', y' \in X$, and $x = x'|A$, $y = y'|B$. The \emph{$h$-almost-union} of $x$ and $y$ is $x \overset{h}\cup y = z \in \Sigma^{A \overset{3r}{\cup} B}$ defined as follows: Define $t \in I$ by setting  $t_a = 0 \iff d(a, A) \leq d(a, B)$. Define $z_a = h(x', y', t)_a$ for all $a \in A \overset{3r}{\cup} B$.
\end{definition}

Note that this depends on the choice of $r$, but we can always use the minimal possible $r$, or alternatively think of $h$ as ``knowing'' its radius. More importantly, note that a priori this definition depends on the choice of extensions $x, y$ of $x', y'$. We show that it does not depend on this choice, and show that it is a globally valid pattern under the assumption contractibility.




\begin{lemma}
\label{lem:PatternAlmostUnion}
The $h$-almost-union $z = x \overset{h}\cup y$ of the patterns $x \in X|A$ and $y \in X|B$ is well-defined, and $z \in X|A \overset{3r}{\cup} B$. Furthermore, if $A, B, x, y$ are $K$-periodic for a subgroup $K \leq G$, then $z$ is $K$-periodic.
\end{lemma}


\begin{proof}
Let $a \in A \overset{3r}{\cup} B$. Let $t$ be as in the definition. We need to show that $h(x', y', t)_a$ does not depend on the choice of $x'$ and $y'$ for 
$a \in A \overset{3r}{\cup} B$. By symmetry, it suffices to show this for $a \in  A^{\circ 3r}$ and for $a \in A \setminus (A B_{3r} \cap B B_{3r}) = A \setminus B B_{3r}$.

If $a \in A^{\circ 3r}$, then either $aB_r$ does not intersect $B^{\circ 2r}$, so $t|aB_R \equiv 0$ and $h(x', y', t)_a = x_a$; or $aB_r$ intersects $B^{\circ 2r}$, and then $aB_r \subset A \cap B$ meaning $h(x', y', t)$ is computed from contents of $x, y$ using the local rule of $h$. If on the other hand $a \in A \setminus BB_{3r}$, then $t|aB_r \equiv 0$ and again $h(x', y', t)_a = x_a$. This shows that $z$ is well-defined.

The pattern $z$ is globally valid for the trivial reason that $z = h(x', y', t)|A \overset{3r}{\cup} B$, where $x', y' \in X$.

If $A, B, x, y$ are $K$-periodic, then since $t$ is chosen by a shift-invariant rule, $t$ is $K$-periodic. Then $h(x', y', t)$ is $K$-periodic in $A \overset{3r}{\cup} B$, since it does not depend on the (possibly non-$K$-periodic) extensions $x', y'$.
\end{proof}

A similar proof gives a strong version of strong irreducibility.

\begin{definition}
Let $G$ be a group, and $K \leq G$. A subshift $X \subset \Sigma^G$ is \emph{$K$-respecting strongly irreducible} if there exists $R$ such that for all subgroups $K \leq G$, and any $K$-periodic $A, B \subset G$ such that $d(A, B) = \min_{a \in A, b \in B} d(a, b) \geq R$, and for any $K$-periodic $x, y \in X$, there exists a $K$-periodic $z \in X$ such that $z|A = x|A, z|B = y|B$. We say $X$ is \emph{period-respecting strongly irreducible} if for all subgroups $K \leq G$, $X$ is $K$-respecting strongly irreducible.
\end{definition}

\begin{lemma}
\label{lem:ContractibleSI}
Let $X \subset \Sigma^G$ be contractible, and suppose the gluing morphism $h : I \times X \times X \to X$ has radius $r$. Then $X$ is period-respecting strongly irreducible with radius $2r+1$.
\end{lemma}

\begin{proof}
Let $K \leq G$. Suppose $A, B \subset G$ are $K$-periodic and such that $d(A, B) \geq R$, and $x, y \in X$ are $K$-periodic. Define $t \in I$ by $t_a = 0 \iff aB_r \cap A \neq \emptyset$. Let $z = h(t, x, y)$. An easy calculation shows that $z|A = x, z|B = y$. Since the four-tuple $(A, B, x, y)$ is $K$-periodic and $t$ is produced by a local rule, we have that $t$ is $K$-periodic. Since $h$ is shift-commuting, it follows that $z$ is $K$-periodic.
\end{proof}

(In fact, $x, y$ only have to be periodic inside the respective domains $A, B$  being glued, but we have opted for the slightly simpler statement.)

\subsection{Contractible implies DPP}

\begin{lemma}
\label{lem:LocalGlobalBlotches}
Suppose $Y \subset A^G$ is contractible and the gluing morphism has radius $r$. Let $x \in A^D$ be a pattern, for some $D \subset G$. If every $2r$-component of $x$ is globally valid, then $x$ is globally valid.
\end{lemma}

\begin{proof}
By Lemma~\ref{lem:ContractibleSI}, the subshift is strongly irreducible with radius $2r+1$. The claim follows by applying strong irreducibility to the components, one at a time, and taking a limit point of the resulting sequence of configurations.
\end{proof}

\begin{theorem}
\label{thm:ContractibleDPP}
Let $Y \subset B^G$ be a contractible subshift on an infinite finitely generated residually finite group $G$. If $G$ has finite periodic asymptotic dimension, then $Y$ has dense periodic points.
\end{theorem}

\begin{proof}
Suppose $G$ has periodic asymptotic dimension at most $d$, where we may assume $d \geq 2$. Contractible subshifts are strongly irreducible by Lemma~\ref{lem:ContractibleSI}. By \cite{CeCo12}, a strongly irreducible subshift containing a point with finite-index stabilizer has dense periodic points. Thus, it suffices to show that there is a periodic point. (Of course, the present proof can be easily modified to directly give dense periodic points.)

Let $Y$ have radius $r$ for the gluing morphism. Since $G$ has finite periodic asymptotic dimension, the subshift $X_{d+1, R, m}$ has a periodic point for arbitrarily large $R$ and some $m$. Let $y$ be any such periodic point, for $R \geq 10 r d$. We may assume that the color $0$ is used, by picking the minimal possible $d$.

The most lengthy part of the proof is showing that we can pick ``centers'' for monochromatic $r$-components in $y$ -- not necessarily in a shift-invariant way, but at least in a periodic way. We begin with this.

By assumption, the group is residually finite. Thus we can pass to an arbitrarily sparse finite-index subgroup $K$ that fixes $y$ (without actually changing $y$), where sparse means that the smallest word norm of a non-identity element is large. We may also assume that $K$ is normal. 

We show that if $K$ is sufficiently sparse, we can pick an element from each color class $K$-periodically. Specifically, we claim that there is a $K$-periodic configuration $x \in \{0, 1\}^G$ such that every $R$-component of a color class $y^{-1}(i)$ of $y$ contains exactly one element $a \in G$ with $x_a = 1$.

We can simply pick centers greedily: as long as there is a monochromatic $R$-component $C \Subset G$ without a center in $y$, we pick a center $a \in C$ for it arbitrarily and then make each element of $Ka$ the center of the component it touches. In other words, we set $x|Ka \equiv 1$. Eventually, every component has a center.

We must check that this process is not contradictory. Specifically, we must check that
\begin{itemize}
\item we never pick a second center for $C$ with this procedure (by shift-symmetry this means we also do not pick two centers for any other component that $Ka$ touches); and
\item no component $D$ touching $Ka$ that gets a center this way had a center previously.
\end{itemize}

For the first item, note that $Ka = aK$ by normality, and
suppose we get a second center for $C$. Then $|C \cap aK| \geq 2$. Then $ab, ac \in C$ for some distinct $b, c \in K$ so $1_G \neq b^{-1}c \in C^{-1}C \cap K^{-1}K$. Since $K = K^{-1}K$, we have that $K$ contains a non-$e_G$ element from $C^{-1}C$. Since $C \subset e B_{Rm}$ for some $e \in G$, $C^{-1}C \subset B_{2Rm}$, and thus we can prevent this by simply picking $K$ sparse enough.

For the second item, if we pick center $a$ for $C$ and $Ka \cap D \neq \emptyset$ (meaning this gives a center for some other existing component $D$), then suppose for a contradiction that $D$ already has a center $b$. We have $KD \cap C \neq \emptyset$ so in fact by $K$-periodicity of $y$ we have $D = cC$ for some $c \in K$. Then $b \in D \implies c^{-1}b \in C$ so $C$ already has a center.

We conclude that we have a $K$-periodic configuration $(x, y) \in (\{0,1\} \times \{0,1,\ldots,d\})^G$ where $x$ marks centers for all $R$-components of $y \in X_{d, r, m}$.


For each color $0 \leq c \leq d$, let $C_c = y^{-1}(c)$, the set of cells $a \in G$ such that $y_a = c$. For each $R$-component $C$ of $C_c$, pick a globally valid pattern $P$ of shape $CB_{4rd}$, whose choice only depends on the shape of the component $C$ (with respect to the unique $1$ in $x|C$). This gives a pattern $x^c$ with support $C_c B_{4rd}$. By Lemma~\ref{lem:LocalGlobalBlotches}, this pattern is globally valid, since the distance between $CB_{4rd}$ and $C'B_{4rd}$ is certainly greater than $2r + 1$ for distinct $R$-components $C, C'$ of $C_c$. Since the patterns $P$ are chosen only based on the shapes of components, $x^c$ is furthermore $K$-periodic.

By induction, the pattern 
\[ y^i = x^0 \overset{h}\cup x^1 \overset{h}\cup \cdots \overset{h}\cup x^i \]
(where we associate to the left) is $K$-periodic. By Lemma~\ref{lem:PatternAlmostUnion}, its support contains
\[ C_0B_{4rd} \overset{3r}\cup C_1B_{4rd} \overset{3r}\cup \cdots  \overset{3r}\sqcup C_iB_{4rd} \]
which contains
\[ ((\cdots (C_0B_{4rd} \cup C_1B_{4rd})^{\circ 3r} \cup \cdots \cup C_{i-1}B_{4rd})^{\circ 3r} \cup C_iB_{4rd})^{\circ 3r} \]
because the $r$-almost-union of two sets contains the union of their $r$-interiors, and thus the $r$-interior of their union. Since in general $(AB_s \cup BB_t)^{\circ r}$ contains $(A \cup B)B_{\min(s, t) - r}$, the support of $y^i$ in particular contains
\[ \bigcup_{j \leq i} C_jB_{rd} \supset \bigcup_{j \leq i} C_j. \]
In particular setting $i = d$, we get that $y^d$ is supported on all of $G = \bigcup_c C_c$. Thus, we have found a $K$-periodic point for a finite-index subgroup $K$.
\end{proof}

We note that the proof is somewhat easier in the case that $G$ is torsion-free, since then no finite component can have a self-symmetry, and we can
immediately pick centers by a shift-invariant rule without needing to pick a new period $N$ and add artificial center markers. (And in this case we can also remove the explicit assumption of residual finiteness by Lemma~\ref{lem:RFC}.)

\begin{question}
On a free group with $k \geq 2$ generators, does there exists a contractible SFT without periodic points?
\end{question}

We suspect that even strong irreducible SFTs have periodic points (and thus dense periodic points \cite{CeCo12a}) on free groups, but we do not have a proof. 

\begin{remark}
Both strong irreducibility and having periodic points are decidable properties for SFTs on free groups (the latter observation is due to Piantadosi, the former we are simply claiming here), so one direction of this problem may be attacked by brute force. We implemented the algorithm for strong irreducibility partially, and Piantadosi's algorithm completely, and generated random SFTs for a few hours without finding a counterexample.
\end{remark}

\subsection{FEP implies DPP}

We have shown that contractibility and SFT imply FEP, and that contractibility implies dense periodic points. In this section we show that FEP also implies dense periodic points. 

\begin{theorem}
\label{thm:FEPDPP}
If $G$ is a finitely-generated infinite residually finite group with finite periodic asymptotic dimension, then FEP subshifts on it have dense periodic points.
\end{theorem}

\begin{proof}
As in the contractible case, since FEP subshifts are strongly irreducible, it suffices to show that one periodic point exists \cite{CeCo12a}. (Though again, one could easily modify the proof so as to directly get a dense set of periodic points.)

Consider an FEP subshift $X$. Let $R$ be the FEP gap, i.e.\ if a pattern on domain $WB_R$ is locally valid (with respect to some fixed finite set of forbidden patterns), then its restriction to $W$ is globally valid.

Pick a generating set $S$ for $G$. Pick some periodic point $y \in X_{d, R, m}$ where $R \geq 10 r d$. For each color $0 \leq c \leq d$, let $C_c$ be the set of cells $a \in G$ such that $y_a = c$. As in the previous proof, we may assume that we also have a periodic binary configuration $x$ that marks a center for each $R$-component of each $C_c$, and that both $y$ and $x$ are $K$-periodic for some finite-index subgroup $K$. 

Note that $3rd < R$, so $r$-components of $C_0 B_{R d}$ correspond (in a bounded-to-one fashion) to $R$-components of $C_0$, and thus are bounded. On each $r$-component of $C_0 B_{r d}$, we position some pattern $P$, where $P$ only depends on the shape of this component -- as in the previous proof, this is possible since the $R$-components have marked centers. 

The pattern of shape $C_0B_{rd}$ obtained is locally valid. By the FEP assumption, the subpattern on $C_0B_{r(d-1)} \subset (C_0B_{rd})^{\circ r}$ is globally valid, and of course it is $K$-periodic. Next, consider $C_1 B_{r d} \setminus C_0 B_{r (d - 1)}$. Again the $r$-components are bounded (they are bounded even before removing the already used cells $C_0 B_{r (d - 1)}$).

For each such $r$-component, we again pick a pattern, only depending on the shape of the component, and the nearby contents of the partial configuration (i.e.\ pattern) we have from the $0$th step. This gives a periodic configuration, which is locally valid. As in the previous step, we shrink it by $r$, to get something globally valid, and continue by induction.

In the next step, we will have a locally valid pattern on
$C_2 B_{rd} \cup C_1 B_{r(d-1)} \cup C_0 B_{r(d-2)}$ followed by a globally valid pattern on $C_2 B_{r(d-1)} \cup C_1 B_{r(d-2)} \cup C_0 B_{r(d-3)}$, and so on.

Finally after the $d$th step (or $2d+1$ if we count the steps where the configuration is only locally valid) we have a periodic configuration defined on all of $\bigcup_c C_c$ and which is $K$-periodic, because the configuration stays $K$-periodic at all times.
\end{proof}


\section{Strong contractibility}

In this section, we define a strengthening of contractibility, which we show to be strongly linked to SFTness and the main definition of \cite{Me23}.

\begin{definition}
We say a subshift $X$ is \emph{strongly contractible}, if there is a homotopy $h : I_2 \times X \times X \to X$ between the left and right projections $\pi_1, \pi_2 : X \times X \to X$, with the additional property that $h(t, x, x) = x$ for all $x \in X$.
\end{definition}

One may also call a strongly contractible subshift one that is contractible mod(ulo) diagonal. 

In the classical topological setting, (the analog of) strong contractibility is called ``equiconnectedness''.

We give a characterization. In topology, a \emph{strong deformation retract} of $Y$ is $X \subset Y$ such that there is a homotopy $h : \ID_Y \to g$ where $g(Y) = X$ and $h(t, x) = x$ for all $x \in X$ and $t \in [0,1]$. This makes perfect sense for subshifts using our symbolic homotopy notion:

\begin{definition}
A subshift $X$ is a \emph{strong deformation retract} of $Y$ if $X \subset Y$ and there is a homotopy $h : \ID_Y \to g$ where $g(Y) = X$ and $h(t, x) = x$ for all $t \in I$ and $x \in X$.
\end{definition}

It is known that in the topological setting, equiconnectedness is equivalent to $X \times X$ strong deformation retracting to the diagonal $\Delta_X$. The same is true in our case, by the same proof, at least when $X$ is SFT.

\begin{lemma}
An SFT $X \subset A^G$ is strongly contractible if and only if $\Delta_X$ is a strong deformation retract of $X \times X$
\end{lemma}

\begin{proof}
Suppose first that $h : I \times X \times X \to X \times X$ is a strong deformation retraction to the diagonal. Define $g(x,y)=\pi_1(h(\bar 1,x,y))$. Then $\pi_1 : X \times X \to X$ is homotopic to $g$ by $h'(t,x,y) = \pi_1(h(t,x,y))$. Similarly, $\pi_2 : X \times X \to X$ is homotopic to $(x, y) \mapsto \pi_2(h(t,x,y))$.

Of course we have $\pi_1(h(1,x,y))=\pi_2(h(1,x,y))$ so $\pi_1$ and $\pi_2$ are two-step homotopic. Furthermore, the homotopies respect the diagonal because $h$ is a strong deformation retraction. Since $X$ is SFT, $\pi_1$ and $\pi_2$ are homotopic by Lemma~\ref{lem:Transitive}. Furthermore it is easy to see that this construction yields a morphisms that respects the diagonal.

On the other hand suppose $\pi_1, \pi_2$ are homotopic mod diagonal by $h : I \times X \times X \to X$. A strong deformation retraction to the diagonal is given by the formula $h'(t,x,y)=(h(t,x,y), y)$.
\end{proof}

Equiconnectedness for topological spaces is not equivalent to contractibility. For example, equiconnectedness implies that for any base point $x_0$, there is a base point preserving homotopy between the space and the constant map at $x_0$. This is called \emph{strong deformation retract contractibility}. It is well-known that there are spaces that are contractible, but one cannot preserve a particular base point. The usual example is the comb space
\[ ((\{1/n \;|\; n \in \N\} \cup \{0\}) \times [0,1]) \cup ([0,1] \times \{0\}) \subset \R^2\]
with the induced topology. There is clearly no base point preserving homotopy between the identity map and the constant map at $(0, 1)$.

Similarly, for subshifts, strong contractibility is not equivalent to contractibility. Indeed, plenty of contractible subshifts are not of finite type, but we show below that on a large family of groups, strong contractibility implies finite type. (On the groups $\Z^d$, we will show the converse as well.)

\begin{definition}
\label{def:GrowthProperty}
Consider a group $G$ with generating set $S$. Let $r, R \in \N$, and consider the smallest family of subsets containing all balls of radius $R$, which is closed under the $r$-almost-union operation and under taking subsets. If for all $r$, there exists $R$ such that left translates $aB_R$, $a \in G$, generate all larger balls in this sense, then we say $G$ has the \emph{patching property}.
\end{definition}


The $r$-almost-union operation is monotone decreasing in $r$, so the choice of metric used to define the property does not matter. The operation is not monotone in $A, B$, but since we explicitly close the family under subsets, we could replace the balls $aB_R$ with left translates of any family of sets $S_i \Subset G$ such that every $b \in G$ is eventually contained in every $S_i$. In particular, the choice of metric does not matter in the definition.


\begin{theorem}
\label{thm:EqContContSFT}
If $G$ has the patching property, then every $G$-subshift $X$ which is strongly contractible is a contractible SFT.
\end{theorem}

\begin{proof}
Contractibility is obvious (even without the patching property).

Suppose $X$ is strongly contractible. Let $h : I_2 \times X \times X \to X$ be the homotopy in the definition of strong contractibility, say it has radius $r/3$. Let $R$ be as in the definition of the patching property for $G$, with parameter $r$, and let $X'$ be the SFT approximation of $X$ where balls $B_R$ of radius $R$ are required to be globally valid for $X$. 

Consider $x \in X'$. Let $\mathcal{C}$ be the family of sets $A$ such that $x|A$ is globally valid in $X$. By definition, the balls of radius $R$ belong to $\mathcal{C}$, and of course this set is closed under taking subsets. It now suffices to show closure under almost-unions, as then every finite pattern in $x$ is globally valid for $X$, meaning $x \in X$.

So consider two $A, B$ such that $x|A$ and $x|B$ are globally valid, say $y, z \in X$ satisfy $y|A = x|A$, $z|B = x|B$. Define $t \in \{0,1\}^G$ by $t_a = 0 \iff d_G'(a, A) \leq d_G'(a, B)$. Here $d_G'(a, C)$ is defined as the distance to $C$ when this distance is positive, and otherwise to be $-n$ where $n$ is maximal such that $aB_n \subset C$.

Then $w = h(t, y, z) \in X$, and it suffices to show that $w$ agrees with $x$ in the set
\[ A^{\circ r} \cup B^{\circ r} \cup (A \cup B) \setminus (B_r(A) \cap B_r(B)), \]
as then this set must be globally valid, thus in $\mathcal{C}$.

Suppose first that $a \in A^{\circ r}$. Then $d_G'(a, A) \leq -r$. If $d_G'(a, B) > -r/3$, then by the $1$-Lipschitz property of $d_G'$, $t|aB_{r/3} \equiv 0$. Thus $w_a = y_a = x_a$ since $h$ is a homotopy. If $d_G'(a, B) \leq -r/3$, then $aB_{r/3} \subset A \cap B$, thus $y|aB_{r/3} = z|aB_{r/3}$, and thus $w_a = x_a$ since $h$ is respects the diagonal. The case $a \in B^{\circ r}$ is symmetric. 

Suppose then that $a \in (A \cup B) \setminus (B_r(A) \cap B_r(B))$. Suppose $a \in A$. Then $d_G'(a, A) \leq 0$ and $d_G'(a, B) > r$, so the first case of the previous argument applies, i.e.\ $w_a = y_a = x_a$. The case $a \in B$ is symmetric.
\end{proof}

\begin{corollary}
Every strongly contractible subshift is SFT on every group with finite asymptotic dimension or subexponential growth. 
\end{corollary}

\begin{proof}
By Lemma~\ref{lem:FiniteDimensionPatching} and Lemma~\ref{lem:SubexponentialPatching}, these properties imply the patching property.
\end{proof}

See Section~\ref{sec:Patching} for some examples of groups covered by this corollary.

On the other hand, on every group, a contractible SFT is strongly contractible as soon as it has a fixed point.

\begin{theorem}
\label{thm:ContractibleFPEqConn}
Every contractible SFT with a fixed point is strongly contractible.
\end{theorem}

\begin{proof}
A contractible SFT $X$ with a fixed point is a retract of a full shift $A^G$, suppose $X \subset A^G$ and let $r : A^G \to X$ be the retraction. Let
\[ c(t, x, y)_a = \left\{\begin{array}{ll}
x_a & \mbox{if } t_a = 0 \\
y_a & \mbox{otherwise.}
\end{array}\right. \]
be the naive homotopy between the left and right projections from $A^G \times A^G$.

Define
\[ h(t, x, y) = r(c(t, x, y)). \]
We have $r(A^G) \subset X$, and for $x, y \in X$ we have
\[ h(\bar 0, x, y) = r(c(\bar 0, x, y)) = r(x) = x, \]
\[ h(\bar 1, x, y) = r(c(\bar 1, x, y)) = r(y) = y, \]
and
\[ h(t, x, x) = r(c(t, x, x)) = r(x) = x \]
for all $t \in I_2$.
\end{proof}


\begin{question}
\label{q:ContractibleSFTEqCon}
Is there a contractible SFT that is not strongly contractible?
\end{question}

We see in Section~\ref{sec:MEP} that for abelian groups, the two properties coincide, so for these groups the answer to the question is ``no''.

\section{Comparison with other dynamical properties}
\label{sec:Comparisons}

\subsection{Contractibility and finite type}

As we have seen, homotopy theory works somewhat more nicely for subshifts of finite type than general ones, in that homotopies are composable in the former case. However, contractible subshifts are not necessarily finite type, and we have shown that density of periodic points holds without this additional assumption (on a large class of groups).

We give here a special property that guarantees contractibility of a subshifts, which allows constructing many non-SFT examples.

\begin{definition}
Let $A \ni 0$ and let $X \subset A^G$ be any subshift defined by forbidden patterns that do not contain the symbol $0$, and whose domains are connected with respect to some generating set $S$ of $G$. Then we say $X$ has property Z0.
\end{definition}

This implies the existence of a safe symbol, and also the property of zero gluing studied in \cite{SaTo21}.

\begin{proposition}
Every subshift with property Z0 is contractible.
\end{proposition}

\begin{proof}
We may use the contraction homotopy
\[ h(t, x, y)_a = \begin{cases}
x_a & \mbox{if } t|aS = 0^{aS}, \\
y_a & \mbox{if } t|aS = 1^{aS}, \\
0 & \mbox{otherwise.} \\
\end{cases} \]

Clearly $h : I \times X \times X \to A^G$ and $h(\bar 0, x, y) = x, h(\bar 1, x, y) = y$. On the other hand if $h(t, x, y)|N$ is one of the forbidden patterns defining $X$, then it must not involve coordinates where $h$ writes $0$, nor can it only involve coordinates coming from only one of the first two cases. Say $h(t, x, y)_a$ comes from the first item (is copied from $x$) meaning $t|aS = 0^{aS}$ and $h(t, x, y)_b$ comes from the second meaning $t|bS = 1^{bS}$.

Since $N$ is $S$-connected, there is a path $p : [0, k] \to G$ with $p(i) \in N$ for all $i \in [0, k]$ and $p(i)^{-1} p(i+1) \in S$. Then there must be a first $i$ such that $t|p(i+1)S \not\equiv 0$. Since $p(i+1) \in p(i) S$, $t_{p(i+1)} = 0$ and thus $t|p(i+1)S$ is not of the form $c^{p(i+1)S}$ for any $c$ in the alphabet, meaning $h(t, x, y)_{p(i+1)} = 0$, contradicting the assumption that $h(t, x, y)|N$ is one of the forbidden patterns defining $X$.
\end{proof}

\begin{proposition}
On every infinite f.g.\ group $G$, there exist uncountably many contractible subshifts over the three-symbol alphabet.
\end{proposition}

\begin{proof}
Let $S$ be the generating set. Let $X \subset \{1, 2\}^\Z$ be any subshift closed under reversal. We define a subshift $Z_X \subset \{0,1,2\}^G$ as a function of $X$. For $z \in A^G$, we have $z \in Z$ if and only if the following holds: for every injective $S$-path $p : I \to G$ where $I \subset \Z$ is a finite interval, if $z_{p(i)} \in \{1, 2\}$ for all $i \in I$, then the $I$-configuration $i \mapsto z_{p(i)}$ is in the language of $X$.

This always gives a $G$-subshift with property Z0 (thus a contractible subshift): the forbidden patterns are injective $S$-paths containing words over alphabet $\{1, 2\}$ not in the language of $X$, in particular these patterns do not contain the symbol $0$.

Using Lemma~\ref{lem:Geodesics} we can find an $S$-geodesic $p : \Z \to G$. If $x \in A^\Z$, then define a configuration $z^x \in A^G$ by setting $z^x_{p(i)} = x_i$ for all $i \in \Z$, and $z^x_a = 0$ for all $a \notin p(\Z)$. Note that $x \mapsto z^x$ is injective.

If $x \notin X$, then clearly $z^x \notin Z_X$ by considering the path $p$. If $X$ is closed under reversal, then $z^x \in Z_X$ for all $x \in X$. Namely, every finite $S$-path in the non-$0$ coordinates of $z^x$ must traverse a finite subpath of $p$ (possibly reversed), and thus sees a word of $X$ or its reversal. Since $X$ is closed under reversal, these are words of $X$.

Thus if $X, X'$ are distinct $\Z$-subshifts closed under reversal, then up to symmetry we have $x \in X \setminus X'$, and $z^x \in Z_X \setminus Z_{X'}$, showing that $Z_X \neq Z_{X'}$. It is easy to show there are uncountably many $\Z$-subshifts that are closed under reversal, concluding the proof.
\end{proof}

The two-symbol alphabet suffices, even using the Z0 property, but the proof is slightly trickier, and we leave this to the interested reader.

\begin{proposition}
On every infinite f.g.\ group $G$, there exists a proper sofic contractible full shift factor.
\end{proposition}

\begin{proof}
Here we use a generating set $S \not\ni 1_G$. For the full shift, we use the alphabet $(S' \times \{0,1\}) \cup \{0,1,\#\}$ where $S' = \{(s, s') \in S^2 \;|\; s \neq s'\}$ meaning each cell contains a bit or $\#$ (which is in fact going to be the zero symbol), and when it contains a bit it can also contain two distinct arrows from the symmetric generating set $S \not\ni 1_G$. The factor map is into the same alphabet. We refer to the $(s, s') \in S^2$ as the \emph{forward} and \emph{backward} arrows, respectively.

The factor map will only modify the bits. Its behavior is as follows: first, we identify \emph{active edges}, namely pairs $(a, as) \in G$ such that the forward edge of $x_a$ is $s$, and the backward edge of $x_{as}$ is $s^{-1}$. The active edges form paths and cycles in the configuration.


Next, we identify \emph{first cells}. These are those cells where a path of positive length through active edges begins; and 
\emph{last cells} where such a path ends. 

In the first cell of a path, we simply write down the current bit. In the last cell, if the forward arrow points to a cell with a bit $b$ and without arrows, and the backward arrow points to a cell with bit $b'$, we write $b \oplus b'$, 
where $\oplus$ is addition modulo $2$. In other active cells -- and also in the last cell if it points to $\#$ or a bit-with-arrows that does not continue the path -- we add the bit in the present cell to the bit in the cell that points to the present cell.

It is easy to check that in the image, on every path of active cells, if the last cell points to a bit without arrows then the sum mod $2$ of the bits on the path, and the bit the last cell points to, is $0$. If the last cell points to $\#$ or a bit-with-arrows, there is no restriction.

The sum constraint on paths where the last cell points to a bit without arrows (which does not continue the path) is in fact the exact subshift constraint defining the image. Thus, the subshift is of type Z0 (with $\#$ the $0$ symbol), so it is contractible. By considering long geodesic segments (obtained from Lemma~\ref{lem:Geodesics}) connected to paths using arrows (and writing $\#$ elsewhere), we can easily find pseudo-orbits that are not actual orbits, so the image is proper sofic.
\end{proof}

\subsection{Contractible SFTs have FEP}
\label{sec:FEP}

In this section, we show that contractibility and finite type together imply the finite extension property. We will see later in Theorem~\ref{thm:NonContractibleFullShiftFactor} that FEP does not imply contractibility.

Recall that once we have fixed a window $W$ for an SFT, \emph{locally valid} for a pattern means there are no forbidden patterns in it (of shape $W$), and \emph{globally valid} means it extends to a full configuration in the SFT.

\begin{theorem}
\label{thm:ContractibleSFThasFEP}
Every contractible SFT has the finite extension property.
\end{theorem}

Note that there is no standard definition of finite extension property for subshifts that are not of finite type, so in that sense the SFT assumption cannot be dropped. 

\begin{proof}
Let $X$ be of finite type, and fix a symmetric window $W$ for it with $1_G \in W$, with forbidden patterns $\mathcal{F} \subset \Sigma^W$. Let $X$ be contractible with homotopy $h : I \times X \times X \to X$ with symmetric neighborhood $R$. Let $\tilde h : I \times \Sigma^G \times \Sigma^G \to \Sigma^G$ be a natural relaxation.

\begin{claim}
Under the above assumptions, for any $x, y \in A^G$, if $\tilde h(t, x, y)|aW \in \mathcal{F}$, then one of the following holds, where $K = aWR$:
\begin{enumerate}
\item $t|K \equiv 0$ and $x|aW$ is a forbidden pattern of $X$;
\item $t|K \equiv 1$ and $y|aW$ is a forbidden pattern of $X$; or
\item $t|K \notin \{0^K, 1^K\}$, and either $x|K$ or $y|K$ is not globally valid in $X$.
\end{enumerate}
\end{claim}

\begin{claimproof}
Suppose $\tilde h(t, x, y)|aW$ is a forbidden pattern. If $t|K \equiv 0$, then $\tilde h$ just copies the left argument, so this forbidden pattern is $x|aW$ and we are in the first item. If $t|K \equiv 1$, we are similarly in the second item.

If $t$ is not constant, we must show $x|K$ or $y|K$ is not globally valid in $X$. Thus, assume to the contrary that $t$ is not constant and $x|K, y|K$ are both globally valid. Let $x', y' \in X$ satisfy $x'|K = x|K, y'|K = y|K$. Then $h(t, x', y')|aW = \hat h(t, x', y')|aW = \hat h(t, x, y)|aW$ is a forbidden pattern, in particular $h(t, x', y') \notin X$, contradicting the assumption about the codomain of $h$.
\end{claimproof}

We now show that under the above assumptions, $X$ indeed has the FEP. For the shape of the forbidden patterns of the FEP we will use the window $WR$ and as forbidden patterns $\mathcal{F}'$ we pick patterns $p \in A^{WR}$ which are not globally valid in $X$. As the FEP gap we will use $N = RW \cup R^2W^2R$. 



If $X$ is empty, the FEP trivially holds, so suppose $X$ is nonempty, and let $y \in X$. Let $M \subset G$ be arbitrary, and let $Q \in A^{M N}$ be any pattern where no pattern from $\mathcal{F}'$ appears. Note in particular that there are then no patterns from $\mathcal{F}$ appearing inside the subpattern $Q|MR$: If $aW \subset MR$, then $aWR \subset MRW \subset MN$, thus $Q|aWR$ is globally valid, thus $Q|aW$ is not a forbidden pattern of $X$.

Now, to prove FEP, it now suffices to show that $P = Q|M$ is globally valid in $X$. For this, let $x \in A^G$ be any configuration with $x|MN = Q$. Define $t$ by $t|MR \equiv 0$, $t|(G \setminus MR) \equiv 1$, and consider $\tilde h(t, x, y) = z \in A^G$. We have $z|M = P$ because $t|MR \equiv 0$ and by our choices of $x$ and the extension $\hat h$. Thus it suffices to show that $z$ is in $X$.

Suppose on the contrary that $z \notin X$. Then $z|aW \in \mathcal{F}$ for some $a \in G$. Let $K = aWR$ and consider the cases of the claim above. The second case (which says $t|K \equiv 1$ and $y|aW \in \mathcal{F}$) is impossible, since $y \in X$. The first case (which says $t|K \equiv 0$ and $x|aW \in \mathcal{F}$) is also impossible: $t|K \equiv 0$ means that $K = aWR \subset MR$ so in particular $aW \subset MR$, but then $x|aW = Q|aW$ would contain a pattern from $\mathcal{F}$ in the subset $MR$, which we argued above is impossible.

Now we look at the third case, i.e.\ $t|K \not\equiv 0$, $t|K \not\equiv 1$, and either $x|K$ or $y|K$ is not globally valid in $X$. Of course it must be $x|K$ that is not globally valid. Since $t$ contains a $0$-symbol in $K = aWR$, $aWR \cap MR \neq \emptyset$ so $a \in MR^2W$ (where we use symmetry of $R, W$). But then $Q$ has a globally forbidden pattern in a translate of $WR$ contained in $K = aWR \subset MR^2W^2R \subset MN$. But this contradicts the assumption on $Q$.
\end{proof}

The reason FEP was originally introduced was the following theorem of Brice\~no, McGoff and Pavlov \cite{BrMcPa18}:

\begin{theorem}
Let $X \subset \Z^d$ be a block gluing subshift, and $Y$ an FEP subshift with a fixed point, such that $h(X) > h(Y)$. Then $X$ factors onto $Y$.
\end{theorem}

Here, recall that a $\Z^d$-subshift $X$ is \emph{block gluing} \cite{BoPaSc10} if there exists $r \geq 0$ such that for every $x, y \in X$, if $A$ and $B$ are half-spaces separated by at least $r$, then there is a point $z \in X$ with $z|A= x|A, z|B = y|B$.

The fixed point is only used to guarantee the existence of a block map $f : X \to Y$. The use of FEP is more subtle, and seems to use the fact $\Z^d$ has finite asymptotic dimension. The point of block gluing and $h(X) > h(Y)$ is to guarantee that some patterns of suitable form in $X$ are able to encode patterns of $Y$. This may be considered as simply a natural way to ensure that there exists a map $g : X \to A^{\Z^d}$ to a full shift, whose image contains $Y$.

With this understanding, the following theorem can be seen as a variant of the above theorem. 

\begin{theorem}
\label{thm:OntoContractible}
Suppose $Y \subset A^G$ is a contractible SFT. Let $X \subset B^G$ be an arbitrary subshift. Then the following are equivalent:
\begin{itemize}
\item There is a block map from $f : X \to Y$ and a block map $g : X \to A^G$ such that $g(X) \supset Y$.
\item There is a factor map $f : X \to Y$.
\end{itemize}
\end{theorem}

\begin{proof}
Obviously the second item implies the first.

Suppose now the first item holds, and let $f, g$ be as in the item. Let $N$ be a symmetric window for $Y$, and let $\mathcal{F} \subset A^N$ be the forbidden patterns of $Y$. Let $h : \{0,1\}^G \times Y \times Y \to Y$ be the contraction homotopy, with symmetric neighborhood $M$. Let $\tilde h : \{0,1\}^G \times A^G \times A^G \to A^G$ be a nastural relaxation. 

Given $x \in X$, define $t(x) \in \{0, 1\}^G$ by
\[ \forall a \in G: (t(x)_a = 1 \iff g(x)|aMN^2M \mbox{ is globally valid in } Y). \]

We now claim that
\[ x \mapsto \tilde h(t(x), f(x), g(x)) \]
is a factor map onto $Y$. First, it is surjective onto $Y$ because if $y \in Y$, then $g(x) = y$ for some $x$. Then $g(x)$ does not contain forbidden patterns of $Y$ so $t(x) = 1^G$, meaning $\tilde h(t(x), f(x), g(x)) = \tilde h(1^G, f(x), g(x)) = g(x)$.

We now explain why the image is contained in $Y$. Suppose it is not, say $\tilde h(t(x), f(x), g(x))|aN = aF$ for some $F \in \mathcal{F}$. Then some symbol in $aN$ must come from $g(x)$ (since $f : X \to Y$). Thus $t(x)_b = 1$ for some $b \in aNM$ (as $\tilde h$ would copy from $f(x)$ otherwise). In particular, by the definition of $t$, $g(x)|bMN^2M \sqsupset g(x)|aNM$ is globally valid in $Y$.

In this case, $\tilde h(t(x), f(x), g(x))|aN = h(t(x), f(x), x')|aN \in Y|aN$ where $x' \in Y$ is arbitrary such that $x'|aNM = g(x)|aNM$.
\end{proof}

\subsection{Non-equality of some classes of SFT}

We have seen that retractions of contractible subshifts are contractible. It is tempting to think that a general factor of a contractible SFT already has nice dynamical properties. A contractible SFT is in particular strongly irreducible, so these factors are strongly irreducible.

In this section we give some indication that one cannot say much more than this, at least using the classes discussed in the present paper. Namely we show:

\begin{theorem}
\label{thm:Separations}
On the groups $\Z^d$, $d \geq 2$, we have the chain of inclusions
\[ \mbox{contractible SFTs} \subsetneq \mbox{FEP} \subsetneq \mbox{SFT factors of contractible SFTs } (\subseteq \mbox{SI SFTs}). \]
\end{theorem}

The first inclusion was already shown in Section~\ref{sec:FEP} (and holds for all groups). The properness of this inclusion is shown in Theorem~\ref{thm:NonContractibleFullShiftFactor} below (the FEP subshift is explicitly given as a full shift factor).

The second inclusion is immediate from \cite{BrMcPa18} for FEP subshifts with a fixed point. We give a general proof in Lemma~\ref{lem:FEPIsCFactor} (for groups $\Z^d$ only). The strictness of the second inclusion is more generally shown for all groups which have a subgroup that is infinite, finitely-presented and one-ended (and the contractible SFT is a full shift). We also show (by a different construction) that on any infinite group that is not just-infinite, there is a separation of the first and third class, i.e.\ a contractible SFT with a non-contractible SFT factor.

The last inclusion is obviously true for all groups, and is included for context only. We suspect this last inclusion is strict as well. At least it cannot be proved non-strict with current knowledge, since if it is not strict, then SI $\Z^3$-SFTs have dense periodic points (which is a well-known open problem).


We need some basic closure properties of contractible subshifts.

\begin{lemma}
\label{lem:FreeExtensionContractible}
If $X$ is a contractible $H$-shift then the free extension to a $G$-shift is contractible.
\end{lemma}

\begin{proof}
We simply apply the homotopy of $X$ on each coset separately. We omit the formal proof. 
\end{proof}

\begin{lemma}
If $X, Y$ are contractible then $X \times Y$ is contractible.
\end{lemma}

\begin{proof}
If $h_X: I \times X \times X \to X$ and $h_Y: I \times Y \times Y \to Y$ are contraction homotopies, then $h : I \times (X \times Y) \times (X \times Y) \to (X \times Y)$ defined by 
\[ h(t, (x, y), (x', y')) = (h_X(t, x, x'), h_Y(t, y, y')) \]
is easily verified to be a contraction homotopy.
\end{proof}

We also need the following.

\begin{lemma}
\label{lem:FEPSFTFreeExtension}
If $X$ is an $H$-subshift and $H < G$, then $X$ is SFT (resp.\ FEP) if and only if $X^{G/H}$ is SFT (resp.\ FEP).
\end{lemma}

\begin{proof}
The SFT case is in \cite{Sa18}.

We prove the FEP claims. We start with the forward implication, which is a little easier. Suppose $X$ is an $H$-subshift and $H < G$. Suppose first that $X \subset A^H$ is FEP, with some defining forbidden patterns $\mathcal{F} \subset A^W$ and FEP window $N \subset H$. Then we can use the same forbidden patterns and the same FEP window for $X^{G/H}$: First, these define $X^{G/H}$ as an SFT essentially by definition. For the FEP property, if $Q \subset A^{MN}$ has no forbidden patterns, then in particular $Q|aH \cap MN$ has no forbidden patterns for any $a \in G$, so $Q|aH \cap M$ is globally valid in $H$, and it follows that $Q|M$ is globally valid in $X^{G/H}$.

We now prove the backward implication. Suppose that $X \subset A^H$ is a subshift, and $X^{G/H}$ is FEP. Suppose $X^{G/H}$ is defined by forbidden patterns $\mathcal{F} \subset A^W$ for some $W \Subset G$ and FEP window $N$. Note that in the definition of FEP, we can always assume that the set of forbidden patterns is just the set of all patterns that are not globally valid in $X^{G/H}$. We can also increase the set $W$ arbitrarily (if we correspondingly increase $N$). Thus we may assume $W = \bigcup_{a \in T'} aW'$ for some $W' \subset H$ and $T'$ a finite set of representatives in distinct left cosets of $H$, one of which is $1_G$.

We claim then that forbidden patterns of shape $W'$ define $X$ as an SFT, and that $X$ is FEP with SFT window $W'$ and FEP window $N \cap H$. For the first claim, if $x \in X$ has no forbidden pattern of shape $aW'$ for any $a$, then extend $T'$ to a complete set of left coset representatives, and define $z \in X^{G/H}$ by $z_{ta} = x_a$ for $t \in T, a \in H$. Then since $W$ defines $X^{G/H}$ as an SFT and consists of translates of $W'$, $z \in X^{G/H}$, in particular we conclude $x \in X$.

For the second claim, suppose $q \in A^{MN}$ does not contain patterns of shape $W'$ that are not globally valid in $X$. Let $D = (G \setminus H) \cup MN$. Pick any $x \in X$, and define $z \in A^D$ by $z_{ta} = x_a$ for $t \in T \setminus \{1_G\}, a \in H$; and $z_a = q_a$ for $a \in MN$. Then for any translate $bW$ of $W$ that fits the domain $D$, either $bW$ does not intersect $MN$ and all the $btW'$ contained in $bW$ see only patterns appearing in $x$; or one of the translates $btW'$ inside $bW$ are contained in $MN$, and (by assumption) contain a valid pattern of $X$, and all other translates $btW'$ contained in $bW$ contain patterns from $x$. In any case, the pattern visible in $bW$ will be valid. By the FEP property of $X^{G/H}$ we conclude that in particular the restriction to $M$, which means that $z|M = q|M$ is valid in $X$, proving the FEP property of $X$.
\end{proof}

\begin{lemma}
\label{lem:FEPIsCFactor}
Let $G = \Z^d$. Then every FEP subshift $Y \subset A^G$ on $G$ is a factor of a contractible SFT.
\end{lemma}
 
\begin{proof}
Let $R$ be the FEP gap, i.e.\ if a pattern on domain $M + B_R$ is locally valid (with respect to some fixed finite set of forbidden patterns), then its restriction to $M$ is globally valid.

Let $D$ be given by Lemma~\ref{lem:ContractibleDimensionSFT}, i.e. such that for any $r'$, the SFT $Y_{D, r', m}$ contains a contractible SFT $Z$ for any large enough $m$.

By Lemma~\ref{lem:XYFactors}, for any $r$, if $r'$ is sufficiently large (as a function of $r$) then we have a morphism from $Y_{D, r', m}$ to $X_{D, r, m'}$ for any sufficiently large $m'$ (as a function of $r'$ and $m$). Take $r > 3RD$. Let the image of the morphism be $X' \subset X_{D, r, m'}$. The product of $Y_{D, r', m}$ with the $A^G$ is still a contractible SFT, and we claim it factors onto $Y$. This is analogous to the proof in \cite{BrMcPa18}, so we give only the outline. (The same idea is also used in Theorem~\ref{thm:FEPDPP}.)

We proceed in steps $i = 1, \ldots, D$, and explain how to produce from any $x \in X'$ and a configuration $z \in A^G$ a valid configuration of $Y$, so that every valid configuration may be obtained. At the beginning of step $i$, we have a partially specified configuration, where for all $j < i$, we have already specified valid contents for the $(r - i R)$-neighborhood of each $j$-component (at $i = 1$, this is trivial -- we have specified nothing in the iamge). Furthermore, at all times the partially specified configuration must have a globally valid extension.

In the inductive step, we pick new contents for each $i$-colored component by reading them from $z$, as follows: We note that the $DR$-neighborhoods around the $i$-components in $x$ are disjoint (and separated by distance at least $R$). In the area not yet determined in previous steps, we now pick the contents from $z$ if this does not give any forbidden patterns, or if there are forbidden patterns, we pick them arbitrarily. This is possible because the configuration is globally valid by the inductive assumption, and we can perform the choice in a shift-invariant continuous way since the components are bounded and disjoint.

Now we shrink the determined area by $R$ steps, to ensure that the remaining configuration is globally valid (using the FEP assumption). Note that in the induction step, we allow the determined area to shrink by $R$, so indeed the induction stays valid. After the final step, we have determined the contents of every cell, and if $z \in A^G$ does not contain forbidden patterns, then we have written exactly $z$.
\end{proof}

\begin{theorem}
\label{thm:NonContractibleFactor}
Suppose $G$ infinite f.g.\ and not just-infinite. Then $G$ admits a contractible SFT which admits a non-contractible SFT factor.
\end{theorem}

\begin{proof}
Let $a \in G \setminus \{1_G\}$. Take finite generating set for $G$ containing $a$. Use colors $\{0,1,2\} \times \{0,1,2\}$. Require that on the first track (i.e.\ the first component $\{0,1,2\}^G$ through the isomorphism $(\{0,1,2\} \times \{0,1,2\})^G \cong \{0,1,2\}^G \times \{0,1,2\}^G$) the symbols $0$ and $2$ cannot be next to each other in any direction; and that on the second track you cannot have the same color between $a$-neighbors. Call the resulting SFT $X$. It is a Cartesian product of two SFTs. The first one has safe symbol $1$ so it is contractible. The second is a free extension of a mixing $\Z$-SFT, so it is contractible.

Now project the second track away when the first track contains $0$ or $2$. The image subshift $Y$ over alphabet $\{0,2\} \cup \{(1,0), (1,1), (1,2)\}$ contains precisely those configurations where on the first track $0$ and $2$ cannot be next to each other, and on the second, on top of two $a$-neighboring $1$s we require values to differ (here one should note that the mixing distance of the one-dimensional SFT we are freely extending is $1$). This is clearly an SFT.

Since $G$ is not just-infinite, $a$ can be chosen so that the normal closure $H$ of $a$ is not of finite index in $G$. We claim that in this case, $Y$ is not contractible. Suppose the contrary, and let $h : I \times Y \times Y \to Y$ be the contraction homotopy. Suppose $t \in \{0,1\}^G$ is constant on cosets of $H$, and consider $x = h(t, \bar 0, \bar 2)$. Since $H$ is normal, every shift of $t$ is $a$-periodic, therefore also every shift of $x$ is $a$-periodic.


Now suppose there are infinitely many cosets for $H$. Let the neighborhood of $h$ be $N$. Let $Ha_1, ..., Ha_k$ be the cosets of $H$ intersecting $N$. Let $t_b = 0$ whenever $b \in \bigcup_i Ha_i$, and $t_b = 1$ otherwise. Now $x = h(t, \bar 0, \bar 2)$ of course contains $0$ at the origin. We claim $x_b = 2$ for some $b \in G$. If this is not the case, then $0 \sqsubset t|bN$ for all $b \in G$. Thus, $bN$ intersects one of the finitely many cosets $Ha_i$ for every $b \in G$. In particular, for every $b \in G$, $b \in Ha_ic$ for some $c \in N^{-1}, i = 1,\ldots,k$, meaning $G \subset \bigcup_{i \leq k, c \in N^{-1}} Ha_ic$ and in fact $H$ has finite index.

Since $x$ contains both symbols $0$ and $2$, it must contain some $(1, \alpha)$ as well, by the first rule of the SFT $Y$. Let $b$ be such that $bx_{e_G} = (1, \alpha)$. Then since $bx$ is $H$-periodic, $abx_{e_G} = bx_{a^{-1}} = (1, \alpha)$ as well, contradicting the SFT rule.
\end{proof}

\begin{theorem}
\label{thm:NonContractibleFullShiftFactor}
If $|A|$ is large enough, then $A^{\Z^d}$ admits an SFT factor which is FEP but not contractible.
\end{theorem}

\begin{proof}
Let $X$ be as in the previous theorem, but with the second $\{0,1,2\}$-component of the alphabet replaced by a larger one (and the rule again that on cosets of some copy of $\Z$, we cannot have adjacent equal symbols). As in the previous theorem, we see that $X$ has a non-contractible SFT factor $Y$.

The subshift$X$ has the single-site filling property, therefore it has the FEP. It also has a fixed point, namely $\bar 0$. It follows that every block-gluing subshift with strictly greater entropy than $h(X)$ factors onto $X$. Thus, for large $|A|$, the full shift$A^{\Z^d}$ factors onto $X$, and thus onto $Y$.
\end{proof}

In fact, every non-trivial full shift $A^{\Z^d}$ admits an SFT factor which is FEP but not contractible -- for this it suffices to find non-contractible FEP subshifts with arbitrarily small entropy and with a fixed point. For this, we can reduce the entropy in the previous construction for example by requiring that symbols $0$ and $2$ all belong to large monochromatic areas, and replacing the binary full shift on $1$-symbols by a low-entropy mixing one-dimensional SFT. We leave the details to the interested reader.

We next give an example that shows that SFT factors of full shifts need not be FEP. This gives another construction of SFT factors of full shifts which are not contractible, in light of Theorem~\ref{thm:ContractibleSFThasFEP}. 

A group is \emph{finitely presented} if it is a quotient of a free group by the normal closure of a finite set. An infinite finitely generated group $G$ is \emph{one-ended} if the following holds: Let $S$ be any finite generating set. Then for every $C \Subset G$, the set $G \setminus C$ has only one infinite connected component in the Cayley graph of $G$.

\begin{theorem}
\label{thm:CohomologicalExample}
Suppose $H$ 
contains a group $G$ which is infinite, finitely presented and one-ended. Then any non-trivial full shift $A^H$ admits an SFT factor which is not FEP.
\end{theorem}

\begin{proof}
It suffices to show that if $G$ is infinite, finitely presented and one-ended, then any non-trivial full shift $A^G$ admits an SFT factor which is not FEP. Namely, if $H$ contains such a group, then we can apply said factor cosetwise to get that $(A^G)^{H/G} = A^H$ has factor $Y^{H/G}$ where $Y$ is SFT and not FEP. Then $Y^{H/G}$ is SFT but not FEP by Lemma~\ref{lem:FEPSFTFreeExtension}.

Furthermore, we may assume the alphabet is $A = \{0,1\}$ by projecting away unnecessary symbols. We now construct a factor map from $A^G$. One should think of the image configurations as having binary values on edges of the Cayley graph. To code this into a literal subshift, we can take alphabet $\{0,1\}^S$ for the image subshift, where $S \not\ni 1_G$ is a symmetric generating set for $G$, 
and we will respect the SFT constraint $(y_a)_s = (y_{as})_{s^{-1}}$ for all $s \in S$ and $a \in G$.

We now simply write the differences modulo $2$ on the edges, i.e.\ if $x \in A^G$ then $(f(x)_a)_s = 1 \iff x_a \neq x_{as}$ (note that then $(f(x)_{as})_{s^{-1}} = 1 \iff x_{as} \neq x_{a}$ so indeed the SFT constraint is satisfied). For algebraic purposes it is more useful to write $(f(x)_a)_s = x_{as} - x_a$.

We claim that the image is SFT, namely we check that the SFT rule defining the image is just that the bits on every cycle sum to $0$, and finite presentation then means it suffices to constrain the parity of finitely many sums.

We check this in more detail. Let $p : \Z_k \to G$ be any cycle such that $p(i)^{-1} p(i+1) \in S$ for all $i \in \Z_k$. Then any $f(x)$ satisfies
\[ \sum_i f(x)_{p(i)} = \sum_i (x_{p(i+1)} - x_{p(i)}) = 0 \]
by telescoping.

Let $R \subset S^*$ be a finite set of relations defining $G$. The geometric meaning is that every $w \in R$ represents a cycle in the Cayley graph of $G$ (when we follow successive prefixes of $w$), and every cycle in the Cayley graph can be written as a concatenation of such cycles, started at different positions, and taking into account cancellations when a path literally moves backwards.

We can then define an SFT where cycles corresponding to $w \in R$ must sum to $0$, and it is then clear that in this SFT, over any cycle, we will see sum $0$, by presenting the cycle as a succession of cycles of shapes from $R$. It follows that $f$ has SFT image.



We claim that the image subshift is not FEP. Suppose on the contrary that $N \Subset G$ is an FEP window, with respect to forbidden patterns fitting inside $W \Subset G$.

Take a geodesic path $p : \Z \to G$ with $\forall i: d(p(i), p(i+1)) = 1$ (so in particular $p(i)^{-1} p(i+1) \in S$). Such a path exists by Lemma~\ref{lem:Geodesics}, and we may assume $p(0) = 1_G$. Now take a large $n$ and observe that $p(i) \in B_n$ (the ball of radius $n$ around the origin) if and only if $i \in [-n, n]$. Both $p(-n-1)$ and $p(n+1)$ have infinite connected components in $G \setminus B_n$, thus they correspond to the same end, and we can join them by a finite path in the subgraph of the Cayley graph of $G$ induced by $G \setminus B_n$.

For some $k$, this gives us a path $q : \Z_k \to G$, where $\Z_k$ is the finite cycle of length $k$, with $q(i)^{-1}q(i+1) \in S$ for all $i$, which agrees with $p$ when $i \in [-n, n]$, and stays outside of $B_n$ on other steps $i \in \Z_k$.

Now consider the configuration $x \in A^G$ defined by $x_a = 1$ if for all large enough $i$ we have $d(p(i), a) \leq i$ (this is the indicator function of the Busemann horoball corresponding to $p$). Clearly $x_{p(i)} = 1 \iff i \geq 0$. Thus on the edges corresponding to the path $p(i)$, we have a single $1$-symbol in the image $f(x)$, namely on the edge from $p(-1)$ to $p(0)$.

Observe that if $j$ has large absolute value, then either $x|p(j)N \equiv 0$ or $x|p(j)N \equiv 1$. Namely if $j < 0$ then
\[ d(p(i), b) \geq d(p(i), p(j)) - d(p(j), b) =  i - j - O(1) > i \]
for large $i$ and for large $|j|$ by the inverse triangle inequality ($O(1)$ accounts for the radius of $N$); and similarly if $j > 0$ then 
\[ d(p(i), b) \leq d(p(i), p(j)) + d(p(j), b) = i - j + O(1) < i \]
for large $i$ and for large $|j|$ by the triangle inequality.

It follows that the restriction of $f(x)$ to the tube around the path $p$, i.e.\ $y = f(x)|p(\Z)N$, has bounded support, i.e.\ the support is contained in a ball, whose radius is independent from the choice of $n$.


Now consider the pattern $z \in (A^S)^{q(\Z_k)N}$ where
\[ z|q([-n, n])N = y|q([-n, n])N = y|p([-n, n])N \]
and $(z_a)_s = 0$ for $a \in q(\Z_k)N$, $s \in S$ not specified by this rule. Then (if $n,k$ are large enough) $z$ does not contain any $W$-shaped forbidden patterns: its values are taken from to $y \in f(A^G)$ in the tube $z|q([-n, n])N$ near the origin, and in the connecting path it corresponds to the ($f$-image of the) constant zero configuration.

The pattern $z$ does not correspond to a restriction of any valid configuration of $f(A^G)$, since along the cycle $q(\Z_k)$, the sum of values on the edges is $1$. We have found a locally valid pattern on $q(\Z_k)N$ whose restriction to $q(\Z_k)$ is not globally valid, a contradiction since $f(A^G)$ was supposed to have the FEP for FEP window $N$.
\end{proof}

The proof of the previous theorem is clearly cohomological, and one can replace finite presentability with the finiteness property that after adding finitely many $G$-orbits of $2$-cells corresponding to relations, the $\Z_2$-valued $2$-cohomology becomes trivial, i.e.\ all cocycles are coboundaries (coboundaries corresponding exactly to configurations with an $f$-preimage).

In the one-dimensional case, all the classes mentioned in Theorem~\ref{thm:Separations} coincide.

\begin{proposition}
Let $G = \Z$. Then an SFT is a strongly irreducible if and only if it is a factor of a contractible SFT if and only if it is FEP if and only if it is contractible.
\end{proposition}

\begin{proof}
It suffices to observe that a one-dimensional SFT is mixing if and only if it is strongly irreducible, as we showed in Theorem~\ref{thm:MixingContractible} that mixing is also equivalent to contractibility.
\end{proof}

\begin{example}
There is a sofic $\Z$-shift which is strongly irreducible but not contractible. Namely, take alphabet $\{0,1,2\}$ and forbid words $2bw2$ where $w \in \{0,1\}^*$ and $b \equiv \sum_i w \bmod 2$. \qee
\end{example}

\begin{question}
Which f.g.\ infinite groups $G$ admit a strongly irreducible SFT which is not a factor of a contractible SFT?
\end{question}

We suspect that ``most'' groups admit such SFTs. At present, we do not know even whether the groups $\Z^d$ do, for $d \geq 2$. It would be difficult to prove that such SFTs do not exist for $\Z^d$ with $d \geq 3$: factors of contractible SFTs (even general contractible subshifts) have dense periodic points by Theorem~\ref{thm:ContractibleDPP}, so showing that their SFT factors are SI would imply that SI SFTs have dense periodic points, which is a well-known open problem. It might be, however, that the classes are easy to separate, although we did not find a construction.

\subsection{Stitching configurations together}

In the case when $X$ is an SFT, contractibility allows us to glue a point together from finitely many configurations, assuming that at least one of them is locally correct everywhere. Further, this can be done in a continuous shift-commuting fashion, meaning we in essence have a way to glue ``partially-defined morphisms'' together. 

\begin{lemma}
\label{lem:Stitching}
Let $X \subset A^G$ be a contractible SFT defined by forbidden patterns $\mathcal{P} \subset A^D$. Let $k \in \N$ be arbitrary. Then there exists a block map $f : (A^G)^k \to A^G$ and $m' \in \N$ such that for all $(x_1, \ldots, x_k) \in A^G$, we have $f(x_1, \ldots, x_k) = x_1$ whenever $x_1 \in X$, and in general we have
\[ \forall a \in G: (\exists i: \forall p \in \mathcal{P}: p \not\sqsubset x_i|aB_{\ell + m} \implies \forall p \in \mathcal{P}: p \not\sqsubset f(x_1, \ldots, x_k)|aB_\ell). \]
\end{lemma}

In other words, $f$ copies its first input when it is locally valid in $X$, and whenever at least one of the inputs has locally valid content from $X$ nearby, then $f$ outputs locally valid content in a slightly smaller area.

\begin{proof}
Observe that by shift-commutation of $f$, it suffices to show the existence of $m$ such that
\begin{equation}
\label{eq:kek}
\exists i: \forall p \in \mathcal{P}: p \not\sqsubset x_i|B_m \implies f(x_1, \ldots, x_k)|D \notin \mathcal{P}.
\end{equation}

Let $h : I \times A^G \times A^G \to A^G$ be a natural relaxation of the gluing morphism, with some neighborhood $M$. Consider first the case $k = 2$. We prove the existence of $f$ in the form~\eqref{eq:kek}. Pick $N$ such that $1 \in N$ and $DM \subset aN$ for all $a \in D$. Then pick $K$ such that $K$ contains $DN$. Finally, by Lemma~\ref{lem:LocalGlobal}, we can find $m$ such that $K \subset B_m$ and
\[ \forall p \in \mathcal{P}: p \not\sqsubset x|aB_m \implies x|aK \sqsubset X \]
for any $x \in A^G$ and $a \in G$.

For $(x, y) \in A^G \times A^G$ define $t(x, y) \in I$ by $t(x, y)_a = 0 \iff x|aN \in X|aN$. This is a block map. Then define $f(x, y) = h(t(x, y), x, y)$.

It suffices to show that if $\forall p \in \mathcal{P}: p \not\sqsubset x|B_m$ or $\forall p \in \mathcal{P}: p \not\sqsubset y|B_m$, then $f(x, y)|D \notin \mathcal{P}$

Suppose first that $\forall p \in \mathcal{P}: p \not\sqsubset x|B_m$. Then $x|K \sqsubset X$, so since $DN \subset K$ we have $x|aN \sqsubset X$ for all $a \in D$. Thus, $t(x, y)|D \equiv 0$, therefore $f(x, y)|D = x|D \notin \mathcal{P}$.

Suppose then that $\forall p \in \mathcal{P}: p \not\sqsubset y|B_m$, so again $y|K \sqsubset X$. If $t(x, y)|D \equiv 1$, then $f(x, y)|D = y|D \notin \mathcal{P}$, so we may suppose $0 \sqsubset t(x, y)|D$, i.e.\ $x|aN \in X|aN$ for some $a \in D$. Since $DM \subset aN \cap B_m$, we have $x|DM, y|DM \sqsubset X$ in this case. Thus there exist points $x', y' \in X$ such that $x'|DM = x|DM$, $y'|DM = y|DM$, and we have $h(t(x, y), x, y)|D = h(t(x, y), x', y')|D \notin \mathcal{P}$ (since $h(t(x, y), x', y') \in X$ by the definition of a homotopy).

When $k > 2$, let $f_2 = f$ the morphism above, and define inductively
\[ f_k(x_1, \ldots, x_k) = f_{k-1}(x_1, \ldots, f_2(x_{k-1}, x_k)). \]
Certainly by induction we have $f_k(x_1, \ldots, x_k) = x_1$ whenever $x_1 \in X$. For the second property, we claim that we can use $(k-1)m$ in place of $m$.

Namely, let $x_1', \ldots, x_{k-2}', x_{k-1}' = x_1, \ldots, x_{k-2}, f_2(x_{k-1}, x_k)$. Suppose for some $i$ we have
\[ \forall p \in \mathcal{P}: p \not\sqsubset x_i|gB_{\ell + (k - 1)m}. \]
Then if $i < k-1$, we of course have have
\[ \forall p \in \mathcal{P}: p \not\sqsubset x_i'|gB_{\ell + (k - 2)m}, \]
and if $i \in \{k-1, k\}$, then we have
\[  \forall p \in \mathcal{P}: p \not\sqsubset x_{k-1}'|gB_{\ell + (k - 2)m} \]
which is the inductive condition for $f_{k-1}$.
\end{proof}

\subsection{Map extension property}
\label{sec:MEP}

Let $X, Y$ be subshifts over the same group $G$. Write $X \leq_{\mathrm{per}} Y$ if for all $x \in X$ there exists $y \in Y$ with $\stab(x) \leq \stab(y)$. This is just the obvious necessary periodic point condition for the existence of a morphism from $X$ to $Y$.

Recently \cite{Me23}, Meyerovitch introduced the following notion.

\begin{definition}
\label{def:MEP}
A $G$-subshift $Y$ has the \emph{map extension property} if the following holds:
Let $X$ be any $G$-subshift such that $X \leq_{\mathrm{per}} Y$. Then for any morphism $g : X' \to Y$ with $X' \subset X$ a subshift, there exists a morphism $f : X \to Y$ that extends $g$.
\end{definition}

Note that $X'$ is allowed to be the empty subshift in this definition. The motivation for introducing this class is quite similar to ours. 
In this section, we give some explicit connections between this definition and contractibility. 
First, we isolate a piece of the definition.

\begin{definition}
We say $Y$ has the \emph{map existence property} if for any subshift $X$ with $X \leq_{\mathrm{per}} Y$, there is a morphism $f : X \to Y$.
\end{definition}

We state the main results of this section, and then prove the lemmas needed in the proofs.

\begin{theorem}
\label{thm:FirstItemCharacterization}
On any finitely generated infinite group $G$, the following are equivalent for a subshift $X$:
\begin{enumerate}
    \item $X$ has the map extension property (Definition~\ref{def:MEP});
    \item $X$ is a contractible SFT with the map existence property;
    \item $X$ is a strongly contractible SFT with the map existence property.
\end{enumerate}
\end{theorem}

\begin{proof}
(1 $\implies$ 3): Suppose $X$ satisfies Definition~\ref{def:MEP}. Then by Lemma~\ref{lem:WMEPImpliesEqCon} below it is a strongly contractible SFT. A subshift with the map extension property has the map existence property because we may take $X' = \emptyset$ in the definition.

(3 $\implies$ 2): A strongly contractible subshift is clearly contractible.

(2 $\implies$ 1): A contractible SFT with the map existence property satisfies the first item of Definition~\ref{def:MEP} by Lemma~\ref{lem:ContSFTHasMEP}.
\end{proof}

\begin{theorem}
\label{thm:ZdFirstItemCharacterization}
On the group $\Z^d$ with $d \geq 2$, the following are equivalent for a subshift $X$:
\begin{enumerate}
\item $X$ has the map extension property (Definition~\ref{def:MEP});
\item $X$ is a contractible SFT;
\item $X$ is strongly contractible.
\end{enumerate}
\end{theorem}

\begin{proof}
The group $\Z^d$ has the patching property, so being strongly contractible is equivalent to being strongly contractible and SFT. By Lemma~\ref{lem:AbelianMExiP}, on the groups $\Z^d$ every contractible SFT (thus also every strongly contractible SFT) has the map existence property. Thus, the items are equivalent to the items of the previous theorem for the groups $\Z^d$.
\end{proof}

\begin{lemma}
\label{lem:WMEPImpliesEqCon}
Every subshift satisfying the map extension property is a strongly contractible SFT. 
\end{lemma}

\begin{proof}
Assume $Y$ has the map extension property. Proposition 4.24 in \cite{Me23} shows $Y$ is of finite type. We show strong contractibility. Take $X = I \times Y \times Y$ and observe that $\pi_2 : X \to Y$, so certainly $X \leq_{\mathrm{per}} Y$.

Let 
\[ X' = (\{\bar 0, \bar 1\} \times Y \times Y) \cup (I \times \Delta_Y) \]
and define $g : X' \to Y$ by the formulas $g(\bar 0, x, y) = x, g(\bar 1, x, y) = y, g(t, x, x) = x$. Then the extension $h = f : X \to Y$ provided by the definition of the map extension property proves strong contractibility.
\end{proof}

\begin{lemma}
\label{lem:ContSFTHasMEP}
Every contractible SFT with the map existence property satisfies the map extension property.
\end{lemma}

\begin{proof}
Assume $Y$ is a contractible SFT and suppose $X \leq_{\mathrm{per}} Y$. We may assume $X \subset A^G$ by picking a common alphabet.

Let $X' \subset X$ be any subshift, and let $g_1 : X' \to Y$ be any morphism. We must show that $g_1$ extends to a morphism $g : X \to Y$. Let $\hat g_1 : A^G \to A^G$ be any extension of $g_1$.

 Let $g_2 : X \to Y$ be some morphism, guaranteed by $X \leq_{\mathrm{per}} Y$ and the map existence property. Let $\hat g_2 : A^G \to A^G$ be any extension of $g_2$.

Let $f : A^G \times A^G \to A^G$ be the stitching map from Lemma~\ref{lem:Stitching}. Define $g(x) = f(\hat g_1(x), \hat g_2(x))$. If $x \in X'$, then $\hat g_1(x) \in Y$ so $g(x) = \hat g_1(x)$ by the first defining property of $f$. Since $g_2(x) \in Y$ whenever $x \in X$, by the second defining property of $f$ we always have $g(x) \in Y$.
\end{proof}

We now show that contractible SFTs on abelian groups have the map existence property. 

We need the usual marker lemma \cite{Kr82}, which states that if a configuration locally does not have a period, then there is a shift-invariant way of dropping ``markers'' on it so that the configuration of markers does not have a period. The following version is proved in \cite{Me23}. We only need the abelian case, and Meyerovitch also proves such version, but we state the general version here and then deduce our own abelian statement.

\begin{lemma}
\label{lem:MarkerLemma}
Let $X \subset A^G$ be a subshift. Let $P \Subset G$ be a finite symmetric set with $1_G \notin P$, and let $V \subset X$ be a clopen set such that $ax \neq x$ for every $a \in P$ and $x \in X$. 
Then there exists a clopen set $C \subset V$ so that
\[ C \cap aC = \emptyset \mbox{ for every } a \in P, \mbox{ and } V \subset C \cup \bigcup_{a \in P} aC \]
\end{lemma}

The clopen set $C$ is referred to as a \emph{marker}. Our viewpoint is that the characteristic morphism $g = \chi_C : X \to \{0, 1\}^G$ defined by $g(x)_a = 1 \iff a^{-1}x \in C$ identifies a set of preferred positions (also called \emph{markers}) in $G$, and these markers are picked in a shift-invariant continuous way.

We state an abelian version in this terminology. If $x \in A^G$ with $G$ abelian, and $a, b \in G$, we say $a$ is a \emph{$b$-period-breaker} if $x_a \neq x_{ab}$.

\begin{lemma}
\label{lem:MarkerLemmaAb}
Let $G$ be an abelian group, let $X \subset A^G$ be a subshift, and let $r, R \in \N$ be arbitrary with $r \leq R$. Then there exists a morphism $g : X \to \{0,1\}^G$ such that for all $x \in X$, $(g(x))^{-1}(1)$ is an $r$-packing (in $G$); and if $x$ has the property that for every $\vec v \in B_r \setminus \{\vec 0\}$, there is a $\vec v$-period-breaker in the set $B_R$, then $(g(x))^{-1}(1)$ is an
$r$-covering of $B_R$.
\end{lemma}

In the statement, $(g(x))^{-1}(1)$ is the set of coordinates in $g(x)$ containing the symbol $1$.

\begin{proof}
Pick $P = B_r \setminus \{\vec 0\}$ and let $V$ be the set of configurations that contain period breakers for all $\vec v \in P$ in the ball $B_{2R}$. Apply the previous lemma and set $g = \chi_C$ where $C$ is the clopen set in the conclusion of the lemma.

The property
\[ C \cap aC = \emptyset \mbox{ for every } a \in P \]
means that $(g(x))^{-1}(1)$ is always an $r$-packing. Now consider the property
\[ V \subset C \cup \bigcup_{a \in P} aC. \]
If in $x$, we can break all $P$-periods in the $R$-ball, then we can break all of them in the $2R$-ball around any $a \in B_R$. This means $ax \in V$ for all $a \in B_R$, and then the condition means that for any $a \in B_R$, $g(x)$ has a $1$ somewhere in $aB_r$.
\end{proof}


Note also that one can alternatively take $V$ to be the set of configurations that see period-breakers in $B_R$, and then use the trick from the proof of  Lemma~\ref{lem:ContractibleDimensionSFT} to spread the markers to the surrounding $R$-ball while preserving the $r$-packing property.


The following notation is nonstandard, but we find it convenient.

\begin{definition}
Let $X \subset A^G$ be a subshift, and $H \unlhd G$ a subgroup. Then write $X/H$ for the $G/H$-system with points $\Fix_X(H)$ and action $aH \cdot x = ax$.
\end{definition}

Often we omit $X$ in the notation $\Fix_X(H)$ when it is clear from context.

\begin{lemma}
\label{lem:QuotientContractibleSFT}
Let $X \subset A^G$ be a contractible SFT. Let $H \unlhd G$ be a normal subgroup. Then $X/H$ is topologically conjugate to a contractible SFT.
\end{lemma}

\begin{proof}
Let $X$ be defined by forbidden patterns $\mathcal{F} \subset A^W$ for some $W \Subset G$. If $p \in \mathcal{F}$ satisfies $p_a = p_{a'}$ whenever $aH = a'H$, then we have a pattern $p/H$ on $G/H$ defined by $(p/H)_{aH} = s \iff p_a = s$. Write $\mathcal{F}/H$ for the set of patterns obtained this way.

Now consider the SFT $X'$ on the group $G/H$ defined by the patterns $\mathcal{F}/H$. We claim that $X'$ is conjugate to $X/H$.

For this, suppose first $x \in X'$. Define $y \in A^G$ by $y_a = x_{aH}$. Since $H$ is normal, we have $y \in \Fix(H)$. Note that
\[ (aH \cdot y)_b = ay_b = y_{a^{-1}b} = x_{a^{-1}bH} = ax_{bH} = (aH \cdot x)_{bH} \]
so we have defined a continuous map intertwining the actions of $G/H$ (note that in the last equality, $ax$ refers to the $G$-action on $X$, and $(aH \cdot x)$ to the $G/H$-action on $\Fix(H)$).

We claim that $y \in X$. Suppose not, say $cy_{W} = p \in \mathcal{F}$ for some $c \in G$. Since $by = y$ for $b \in H$, a short calculation shows $p_a = p_{a'}$ whenever $aH = a'H$ and $a, a' \in W$. Then a short calculation shows that $(cH \cdot x)|WH = p/H \in \mathcal{F}/H$, a contradiction.

We can similarly construct a map in the other direction, and verify that these maps commute.

To see contractibility, let $h : I_G \times X \times X \to X$ be the contraction homotopy, and define $h' : I_{G/H} \times X/H \times X/H \to X/H$ using the obvious pullback formula, explicitly $h'(t', x', y')_{bH} = h(t, x, y)_b$ where $t_a = t'_{aH}, x_a = x'_{aH}, y_a = y'_{aH}$ for all $a \in G$. This is well-defined, because by $G$-commutation, the image of $h$ is in $\Fix(H)$ whenever its inputs are.
\end{proof}

\begin{remark}
We note the subtlety that $\Fix(H) \subset X$ is not necessary a subshift of finite type as a subset of $A^G$, i.e.\ the $G$-subshift obtained as the image of $X'$ in the above proof is not an SFT for the full $G$-action, only for the $G/H$-action. For example, let $X$ be the two-symbol full shift on the lamplighter group $\Z_2 \wr \Z$. Then $\Fix(\bigoplus_\Z \Z_2)$ is a free extension of the two-point subshift on $\bigoplus_\Z \Z_2$. Since this latter is not SFT (\cite{Sa18}).
\end{remark}

\begin{lemma}
\label{lem:AbelianMExiP}
Every contractible $G$-SFT with $G$ finitely-generated abelian has the map existence property.
\end{lemma}

\begin{proof}
Suppose $Y \subset B^G$ is a contractible $G$-SFT.

We may inductively suppose that for all proper quotients of $G$, the claim is true, since f.g.\ abelian groups have the ascending chain condition for subgroups. 

Suppose $X \leq_{\mathrm{per}} Y$ for two $G$-subshifts. For $G$ the trivial group, the claim is trivial: every $x \in X$ is simply mapped to the same configuration $y_0 \in Y$.

For the general case, pick a symmetric generating set $S$. Let $r'$ be large. By default, we will apply the previous marker lemma (with $r = r'$) to identify an $r'$-net of markers (in the sense that the set is $r'$-separated and is $r'$-dense other than near the boundary). 


In regions containing such markers, we look at the Voronoi cells of the markers. They give us patterns satisfying the SFT constraint of the subshift $Y_{D, r', m'}$ from Lemma~\ref{lem:XYFactors} for sufficiently large $D$ and $m'$. Note that $D$ can be taken independent from $r'$ as in the proof of Lemma~\ref{lem:ContractibleDimensionSFT}, by a volume argument and using the fact $G$ is doubling as a metric space.

We can then construct points of $X_{D, r, m}$ using the proof of Lemma~\ref{lem:XYFactors}, for some $r$. As $r' \rightarrow \infty$, $r$ can also be taken arbitrarily large. We then follow the proof of Lemma~\ref{lem:FEPIsCFactor} to construct valid contents in $Y$, by proceeding one color at a time. This is possible since $D$ stays bounded while $r$ can be taken arbitrarily large. 

Extending this map arbitrarily, we obtain a morphism $f_{\vec 0} : X \to B^G$ such that when the marker lemma gives us markers sufficiently nearby (meaning there are period-breakers sufficiently nearby, say at distance at most $R/2$), then $f_{\vec 0}$ does not produce forbidden patterns of $Y$ nearby.

In areas where the marker lemma does not apply, i.e.\ large areas with at least one small period, we apply induction. Specifically, we note that for any vector $\vec v$, $Y/\langle \vec v \rangle$ is (by the previous lemma) a contractible SFT on the proper quotient group $G/\langle\vec v\rangle$. Thus, it has the map existence property. In general if $H \lhd G$, and $X \leq_{\mathrm{per}} Y$, then one can verify $X/H \leq_{\mathrm{per}} Y/H$, therefore we have a morphism $f_{\vec v}' : X/\langle\vec v\rangle \to Y/\langle\vec v\rangle$.

From $f_{\vec v}'$ we obtain a morphism
\[ f_{\vec v}'' : \Fix_X(\langle \vec v\rangle) \to \Fix_Y(\langle \vec v\rangle) \]
in an obvious way, and we can then extend this further to a map $f_{\vec v} : X \to B^G$ arbitrarily. We observe that when $x \in X$ is $\vec v$-periodic in a large enough area, then $f_{\vec v}$ produces a valid pattern from $Y$ in a large area.



All in all, we have morphisms $f_{\vec 0}$, and $f_{\vec v}$ for each $\vec v \in P$. We observe that without changing $r'$ (thus with a fixed number $|P| + 1$ of morphisms), we can ensure that for each $\vec u \in \Z^d$, an arbitrarily large ball is covered in the valid part of the image of at least one of the morphisms $f_{\vec v}$, $v \in B_{r'}$. For this, it suffices to increase the number $R$, so as to make the morphism $f_{\vec 0}$ look further for period-breakers.

All in all, we conclude that for any $\ell$, we can ensure that for all $x \in X$, at least one of the morphisms $f_{\vec v}$ produces a valid pattern of $Y$ in the ball of radius $\ell$. We now obtain a morphism from $X$ to $Y$ using Lemma~\ref{lem:Stitching}.
\end{proof}





\section{The new group properties}

In the main text, we abstracted away some group properties that we needed in proofs, namely the finite periodic asymptotic dimension, and the patching property. To our knowledge, these notions have not previously appeared in the literature. Here, we give examples of groups with these properties.

\subsection{Groups with finite periodic asymptotic dimension}
\label{sec:Dimension}

\begin{theorem}
\label{thm:ZdDimension}
The periodic asymptotic dimension of $\Z^d$ is $d$.
\end{theorem}

\begin{proof}
The usual proofs that the asymptotic dimension of $\Z^d$ is at most $d$ already use a periodic grid.

Let us recall this argument. We should prove that for all $r$, there exists $m$ such that the $d$-dimensional subshift $X_{d+1, r, m}$ has a periodic point. 

Let $N$ be large, and for $\vec v \in \Z^d$, denote by $k(\vec v)$ the maximal $k$ such that there are at least $k$ coordinates $i$ such that $d(\vec v_i, N\Z) \leq kr$. Let $C_{k'} = \{\vec v \;|\; k(\vec v) = k'\}$, and observe that $\{C_k \;|\; 0 \leq k \leq d\}$ is a partition of $\Z^d$.

Suppose $\vec u$ and $\vec v$ are in $C_k$ and the $\ell_\infty$ distance between $\vec u$ and $\vec v$ is at most $r$. Let $I$ be the $k$ many coordinates such that $d(\vec u_i, N\Z) \leq kr$ for $i \in I$, and let $J$ be the corresponding set for $\vec v$. Observe that for $i \notin J$, we must have $d(\vec v_i, N\Z) > (k + 1)r$, or $k$ is not maximal for $\vec v$. 

Then we must have $I = J$. Namely, suppose on the contrary that $|I| = |J|$ and the distance is at most $r$, but $I \neq J$. If $i \in I \setminus J$, then we would have
\[ kr <  d(\vec v_i, N\Z) - r \leq d(\vec u_i, N\Z) \leq kr, \]
a contradiction.

We conclude that any path $p : \N \to T_k$ that satisfies $d(p(i), p(i+1)) \leq r$ must satisfy that there is a set $I$ of $k$ coordinates such that $\vec v = p(i)$ has $d(\vec v_i, N\Z) \leq kr$ for $i \in I$ and $d(\vec v_i, N\Z) > (k + 1)r$ for $i \notin I$.

If $N > 10(d+1)r$, it is easy to see that in each individual coordinate, the set of possible values it can reach is finite. Namely if $i \in I$, then we stay in a single hyperplane $U_i = \Z^{i-1} \times [cN - kr, cN + kr] \times \Z^{d-i-2}$, since such hyperplanes are at distance more than $r$ for distinct $c$; and if $i \notin I$, then we stay in a single hyperplane $U_i' = \Z^{i-1} \times [cN + (k+1)r, (c+1)N - (k+1)r] \times \Z^{d-i-2}$, since such hyperplanes are at distance more than $r$ for distinct $c$

This gives $\pdim(\Z^d) \leq d$. For the lower bound $\pdim(\Z^d) \geq d$, we simply observe $\pdim(G) \geq \dim(G) = d$ where the latter is well-known.
\end{proof}

\begin{lemma}
\label{lem:SubgroupPDim}
Let $H \leq G$ be finitely-generated groups. Then $\pdim(H) \leq \pdim(G)$.
\end{lemma}

\begin{proof}
If $x \in X_{G, d, r, m}$, then clearly $x|H \in X_{H, d, r, m}$ (taking the generating set of $G$ to contain the one used for $H$).
\end{proof}

\begin{lemma}
\label{lem:FiniteIndex}
If $H < G$ is of finite index and both groups are finitely-generated, then $\pdim(H) = \pdim(G)$.
\end{lemma}

\begin{proof}
The group $H$ cannot have larger periodic asymptotic dimension than $G$ by the previous lemma. We now show that $G$ does not have larger periodic asymptotic dimension than $H$, which concludes the proof.

Let $S$ be the finite generating set used for $H$, $S' \supset S$ the one for $G$, pick coset representatives $T$ so $G = \bigsqcup_{t \in T} Ht$.
For each $t \in T, s' \in S'$ write $ts' = s''t'$ with $s'' \in H, t' \in T$. Let $\ell$ be the maximal word norm of any $s''$ that appears this way.


Now if $x \in X_{H, d, r, m}$ is periodic with stabilizer $K \leq H$, then we construct a point $y \in \{0, 1, \ldots, d\}^G$ by the formula $y_{hs} = x_h$ for $h \in H, s \in T$. Then the stabilizer of $y$ still contains $K$ which is of finite index in $G$.

We claim that $y \in X_{G,d,r',m'}$ for some $r',m'$, where $r' \rightarrow \infty$ as $r \rightarrow \infty$. Namely, consider a monochromatic $r'$-path in $y$. Then we can construct a \emph{shadow path} in $H$ by projecting $ht$ to $h$ for $t \in T, h \in H$. Consider an $S'$-move from $ht$ to $hts'$ with $s' \in S'$. We have $hts' = hs''t'$ with $s'' \in H$, $t' \in T$ with $s''$ of length at most $\ell$. Thus, a monochromatic $r'$-path with $S'$-moves in $G$ projects to an $r'$-path in $x$ with $S$-moves. So we can pick $r' = \lfloor r/\ell \rfloor$ and $m' = m|T|$.
\end{proof}

A group $K$ is \emph{locally finite} if every finitely-generated subgroup of $K$ is finite.

\begin{proposition}
\label{prop:LocallyFiniteKernel}
Let $1 \to K \to G \to H \to 1$ be an exact sequence, and suppose $K$ is locally finite. Then $\pdim(G) \leq \pdim(H)$.
\end{proposition}

Note that for a split extension, we then have $\pdim(G) = \pdim(H)$ by Lemma~\ref{lem:SubgroupPDim}.


\begin{proof}
For $G$ use a generating set $S$ which projects to the generating set used for $H$. Let $\pi : G \to H$ be the projection with kernel $K$.

For $r \in \N$, let $m$ be such that there is a periodic $x \in X_{H,d,r,m}$ where $d = \pdim(H)+1$. Pull this back to a $G$-configuration through the quotient map by $y_a = x_{\pi(a)}$. Then $y$ obviously has finite orbit. We claim that $y \in X_{G,d,r,m'}$ for some $m'$. For this, observe that while staying within a monochromatic component, the projection will stay in the same component of $x$. It suffices to give a bound on the size of the components in the fiber $F' = \pi^{-1}(F)$ where $F$ is the monochromatic component of $x$ that contains the identity.

Consider a path by right $S$-translations which stays stays in the fiber $F'$. Picking a lift of $F$ and canceling and reintroducing the right coset representative of elements encountered during the path as in the previous lemma, we get a shadow path which stays inside $K$, and moves by finitely many generators $T \Subset K$. Since $K$ is locally finite, there is a bound $m'$ on the components. This bound can be taken to only depend on the diameter of $F$, so we are done.
\end{proof}

\begin{corollary}
The lamplighter group $G = \Z_2 \wr \Z$ has $\pdim(G) = 1$.
\end{corollary}

We now consider more general group extensions.

\begin{theorem}
\label{thm:GroupExtensionDimension}
Let $G$ be a finitely-generated group. Suppose $1 \to K \to G \to H \to 1$ is split exact, and every finite-index subgroup of $K$ has finite orbit under the conjugation action of $H$, then $\pdim(G) \leq (\pdim(K) + 1)(\pdim(H) + 1) - 1$.
\end{theorem}

The assumption about finite-index subgroups is implied by $K$ being finitely-generated, but we give the more general statement to cover the Baumslag-Solitar group $\BS(1, n)$ below.

\begin{proof}
Recall that in the case that a group $K$ is not finitely generated, in the definition of the dimension, by $r$-paths we simply mean paths by jumps with sets $S_r$, where $K = \bigcup_r S_r$.

Let $\pi : G \to H$ be the quotient map with kernel $K$. Use on $G$ a generating set that projects onto the generating set of $H$. We need $\pdim(K) + 1$ colors $D$ for the colorings of $K$, and $\pdim(H) + 1$ colors $D'$ for the colorings of $H$. We will use colors $D \times D'$ for $G$ to get the stated dimension formula.

Let $r$ be given. We construct a periodic point $y$ in $X_{G, |D \times D'|, r, m}$ for some $m$. First, take a periodic point $x \in X_{H,D',r,m'}$ for sufficiently large $m'$. If $\pi : G \to H$ is the projection, set the $D'$-component of $y_a$ to be $x_{\pi(a)}$.

Next, consider a fiber corresponding to a single finite $r$-component $F \Subset H$ of a color class of $x$. Then no matter how we pick the $D$-components in $y$, all monochromatic $r$-components of elements in the fiber $\pi^{-1}(F)$ stay in this fiber (since the distance in $G$ between two distinct fibers is larger than $r$, because of the compatibility of the generating sets of $G$ and $H$). Thus we are effectively left with picking the $D$-colors in individual fibers, so that $r$-components within fibers are bounded, and the global configuration is periodic.

Let us consider the fiber $F$ containing the identity element. As in the proof of Proposition~\ref{prop:LocallyFiniteKernel}, we observe that a path by $r$-jumps in $G$, which stays inside the fiber $\pi^{-1}(F)$, stays close to a path in $K$ by jumps from a finite set. Specifically, let $R$ be a (finite) set of representatives for $\pi^{-1}(F)$, then any path inside $\pi^{-1}(F)$ with $r$-jumps moves by jumps from $kt, k't'$ for $k, k' \in K$ and $t, t' \in R$, and we can again define its shadow path as the sequence of elements in $K$. In the shadow path, we thus move by elements of $RB_rR^{-1} \cap K$.

Let now $r'$ be such that the set $S_{r'}$, used to define dimension of $K$, contains $RB_rR^{-1} \cap K$. There exists $m'$ such that $X_{K, |D|, r', m'}$ has a point $z$ with finite index stabilizer $N_z$. We can then copy the color at $k \in K$ to each $kt$ with $t \in R$, to obtain a $D$-coloring of $\pi^{-1}(F)$. Observe that $K$ acts by translation on patterns with domain $KR$ (since $K \cdot KR = KR$) and copying the color from $K$ to the right cosets $KR$ preserves the stabilizer $N_z$. From this discussion, we conclude that there is a $D$-coloring of $\pi^{-1}(F) \subset G$ which has finite-index stabilizer, and where monochromatic components are bounded.

We can use a similar scheme to pick the $D$-symbols in any individual fiber $\pi^{-1}(F)$ of a monochromatic $r$-component $F \Subset H$ by translating it so that it contains the identity. Note that there are only finitely many sets $F$ that can appear, since these components are connected and bounded. No matter how we do this, it is automatic that the $r$-components are bounded if we use only finitely many types of coloring $m$ (which we can do, since there are finitely many choices of $F$). We are left with making the choices so that the global configuration is periodic.

Since the extension $1 \rightarrow K \rightarrow G \rightarrow H \rightarrow 1$is split, there is indeed a canonical way to perform a coloring, which will lead to periodicity. First, we may assume that the stabilizer $L$ used for $x$ is normal and sparse enough, so that, as in the proof of Lemma~\ref{lem:ContractibleDimensionSFT}, we can pick ``centers'' for the finite monochromatic connected components $F \Subset H$, so that the set of centers is $L$-periodic.

Now the scheme we use to color $\pi^{-1}(F)$ is that we take the center $h \in F$, see $h$ as an element of $G$, translate it to the origin, color the fiber, and move it back.

The canonicality of this procedure implies that the resulting configuration is automatically $L$-periodic. We claim that its $K$-orbit is also finite. To see this, note that we use only finitely many different periodic points $z$ in the construction, and thus we have only finitely many stabilizers $N_z$. Let $N$ be the intersection of these groups. By the assumptions, $N$ has finite orbit under the conjugation action of $H$. Let $L' \leq N$ be the intersection of all these conjugates $L' = \bigcap_{a \in H} a^{-1}Na$. Then $L'$ is a finite-index subgroup that is invariant under conjugation by elements of $H$.

Observe that a translation by $k$ is effectively translation by $a^{-1}ka \in K$ in the fiber $\pi^{-1}(aF)$, in the sense that if the contents of $\pi^{-1}(aF)$ are obtained by taking $z \in D^{\pi^{-1}(F)}$ and translating it to $\pi^{-1}(aF)$ by shifting $z \mapsto az$, then $k az = a (a^{-1}k a) z$ meaning the new contents of $\pi^{-1}(aF)$ corresponds to the $a$-shift of $(a^{-1}k a) z$. For $k \in L'$, $(a^{-1} k a) z = z$ because $a^{-1}ka\in L'$ and each $z$ is $L'$-periodic (which in turn holds because $L'$ contains $N$, which contains $N_z$).

We conclude that the stabilizer of the configuration we have built projects onto a finite index subgroup (namely $L$) of $H$, and contains a finite index subgroup (namely $L'$) of $K$. It is a general group-theoretic fact that such a subgroup has finite index in $G$, thus we have constructed a periodic point in $X_{G, |D \times D'|, r, m''}$ for some $m''$.
\end{proof}

Note that whenever there is a finite number of subgroups of any given finite index, the orbit of any finite-index subgroup must be finite, since the conjugation action of $H$ is by automorphisms of $K$. Furthermore, it is well-known that when $K$ is finitely-generated, then it has finitely many subgroups of any given finite index. Thus, the theorem covers all cases where $K$ is finitely-generated, which leads to the following corollary.

\begin{corollary}
Virtually polycyclic groups have finite periodic asymptotic dimension.
\end{corollary}

\begin{proof}
A strongly polycyclic group is inductively defined as a group extension $G$ with exact sequence $1 \rightarrow P \rightarrow G \rightarrow \Z \rightarrow 1$ where $P$ is itself strongly polycyclic or trivial. Since $\Z$ is a free group, any such extension splits, so the previous lemma and induction shows that they have finite periodic asymptotic dimension.

Now the result follows from Lemma~\ref{lem:FiniteIndex}, since a virtually polycyclic group by definition has a finite-index subgroup which is polycyclic, and it is well-known that a polycyclic group has a strongly polycyclic finite-index subgroup.
\end{proof}

\begin{remark}
As we saw in the proof of Theorem~\ref{thm:GroupExtensionDimension}, the assumption ``every finite-index subgroup of $K$ has finite orbit under the conjugation action of $H$'' equivalently says that every finite-index subgroup contains a finite-index subgroup that is invariant under the action of $H$ (as explained in the proof), equivalently invariant under the conjugation action of $G$. Thus, another way to  state the condition is by saying that among finite-index subgroups of $K$, the ones that are normal in $G$ are cofinal under reverse inclusion.
\end{remark}

As we have already mentioned, it is well-known that finitely-generated groups have finitely many subgroups of a given finite index. In exactly the same way we can prove a slight generalization.

\begin{lemma}
Let $K$ be locally boundedly generated, meaning for some $k \in \N$, every finitely-generated subgroup admits a generating set with $k$ elements. Then $K$ has finitely many subgroups of any given finite index.
\end{lemma}

\begin{proof}
Let $S_i$ be sets of cardinality $k$ such that $S_i \subset \langle S_{i+1}\rangle$ and $K = \bigcup_i \langle S_i \rangle$.

The number of subgroups of index $n$ is at most the number of actions on $\{1, \ldots, n\}$ (because any subgroup gives a transitive action on the set of left cosets). The number of ways $S_i$ can act is at most $n!^k$. If we have more than $n!^k$ distinct actions of $K$, then by the pigeonhole principle, two of them will use the same action of $S_i$ for infinitely many $i$. Since $K = \bigcup_i \langle S_i \rangle$, the actions of $K$ must then be the same.
\end{proof}

\begin{corollary}
The Baumslag-Solitar group $\BS(1, n)$ has periodic asymptotic dimension $2$ or $3$.
\end{corollary}

\begin{proof}
The group $\BS(1, n)$ is a semidirect product $\Z[\frac1n] \rtimes \Z$ (equivalently, a split extension) where $G = \Z[\frac1n] = \bigcup_k \frac1{n^k} \Z$ as an additive subgroup of $\Q$. The quotient $\Z$ of course has dimension $1$.

The group $G$ is locally cyclic, so by the previous lemma it has finitely many subgroups of any given finite index, thus the orbit of any finite-index subgroup is finite under the entire automorphism group.\footnote{Actually, one can show that every finite index subgroup of $G$ is characteristic.}

We now show that $G = \Z[\frac1n]$ has periodic asymptotic dimension at most $1$ (and thus precisely $1$). For this, let $r$ be arbitrary, and $S_r$ the set used to define the dimension of $G$. Up to applying an automorphism of $G$ (the automorphism group of $G$ is isomorphic to $\Z \times \Z_2$ and acts by multiplication by powers of $2$ and by $\pm 1$), we may assume that $S_r$ generates $\Z$. Suppose further that $S_r \subset [-r', r'] \cap \Z$.

Take $p$ a prime number larger than $10r'$ that does not divide $n$, and consider the homomorphism $\pi : a/n^\ell \mapsto a/n^\ell \bmod p$ from $G$ to $\Z_p$, where in $\Z_p$ we interpret division as using the field structure. This is well-defined, since $n$ is nonzero in $\Z_p$, and we have clearly defined a homomorphism. Define $x \in \{0,1\}^{\Z_p}$ by $x_i = 0 \iff i \in \{-r', \ldots, r'\}$. Then the pullback $y \in \{0,1\}^G$ defined by $y_i = y_{\pi(i)}$ has orbit of size $p$, and it belongs to $X_{G, 2, r, m}$ for $m = \max(p - (2r'+1), 2r'+1)$, since with $S_r$-generators we cannot move from one fiber of $\Z_p$ to another along a monochromatic component.

The previous theorem now gives the bound $3$ for the dimension, and the well-known asymptotic dimension $2$ of the Baumslag-Solitar group gives the lower bound.
\end{proof}

Presumably, the periodic asymptotic dimension of the Baumslag-Solitar group is precisely $2$.

Generally, one expects dimension to behave according to $\dim(A \times B) \leq \dim(A) + \dim(B)$, where $A \times B$ is some type of product of $A$ and $B$. It is known that in the case of asymptotic dimension of groups, such a formula holds for group extensions \cite{BeDr06}, and in the case of Lebesgue covering dimension of metric spaces, it holds for the direct product.

Our naive formula suffices for us because in our applications we are only interested in finiteness of the dimension. We do not know in which situations the stronger formula holds. 

\begin{question}
Does periodic asymptotic dimension stay finite in general group extensions? Does the formula $\pdim(G) \leq \pdim(K) + \pdim(H)$ hold for periodic asymptotic dimension for (split or non-split) group extensions?
\end{question}

\begin{conjecture}
Nonabelian free groups have infinite periodic asymptotic dimension.
\end{conjecture}

Note that the usual asymptotic dimension of free groups is $1$. The groups $\BS(m, n)$ with $1 < m \leq n$ contain nonabelian free subgroups, so if the conjecture is true, then they are also infinite-dimensional.

\subsection{Groups with the patching property}
\label{sec:Patching}

We have no examples of groups that do not have the patching property from Definition~\ref{def:GrowthProperty}. However, we give two (incomparable) classes of groups that do have it. 

These cover for example free groups, virtually polycyclic groups, metabelian Baumslag-Solitar groups, the lamplighter group $\Z_2 \wr \Z$, free groups, and the Grigorchuk group.

\begin{lemma}
\label{lem:FiniteDimensionPatching}
Let $G$ be a f.g.\ group with finite asymptotic dimension. Then $G$ has the patching property.
\end{lemma}

Note that here we talk about asymptotic dimension rather than periodic asymptotic dimension, so this covers the free group.

\begin{proof}
Take $C_1 \cup \cdots \cup C_d$ a cover of the group so that each $C_i$ has connected $2r$-components bounded, and where each element $a \in G$ has its $2rd$-ball contained in some $C_i$ (for this, start from slightly more separated sets $C_i$, and then artificially increase their sizes).

Then balls of radius $B_R$, where $R$ comes from the bound on the size of connected components, give the patching property for $r$. First, any finite subset of the union of connected components in a given $C_i$ can be constructed using almost-unions, since $B_r(A) \cap B_r(B) = \emptyset$ for distinct connected components $A, B$, so for distinct components
\[ A \overset{r}{\cup} B \supset (A \cup B) \setminus (A B_r \cap B B_r)) = A \cup B \]

Now let $N \Subset G$ be arbitrary, and build a set $A_i \Subset C_i$ such that $(N \cap C_i^{\circ dr}) \subset A_i^{\circ dr}$. Then
\[ A_1 \overset{r}{\cup} A_2 \overset{r}{\cup} \cdots \overset{r}{\cup} A_d \supset A_1^{\circ rd} \cup A_2^{\circ rd} \cup \cdots \cup A_d^{\circ rd} \]
By the assumption on the $C_i$, each $a \in N$ satisfies $aB_r \subset C_i$ i.e.\ $a \in C_i^{\circ dr}$, for some $i$, thus $a \in A_i^{\circ dr}$ for some $i$ by the choice of the sets $A_i$. We conclude that we have built a superset of $N$.
\end{proof}

\begin{lemma}
\label{lem:SubexponentialPatching}
Let $G$ be a finitely-generated group with subexponential growth. Then $G$ has the patching property.
\end{lemma}

\begin{proof}
Let $r \in \N_+$ be arbitrary. Suppose that there exists $k$ such that $B_{R+1}$ can be covered by at most $2^k$ many left translates of $B_{R - rk}$, for all large enough $R$. Then we claim that $R$ proves the patching property for $r$, i.e.\ translates of $B_R$ generate arbitrarily large sets under $r$-almost-unions.

Namely, let us cover $B_{R + 1}$ with $2^k$ many translates $a_i B_{R-rk}$, $a_i \in G$. Now combining the corresponding translates $a_i B_R$ pairwise in a tournament-tree fashion, we get a common almost-union for them with at most $2^k-1$ almost-unions, so that each individual $a_i B_R$ only partakes in an almost-union on $k$ steps.

The almost-union of $A$ and $B$ contains $A^{\circ r} \cup B^{\circ r}$ and the operation $A \mapsto A^{\circ r}$ is superadditive in the sense that $A^{\circ r} \cup B^{\circ r} \subset (A \cup B)^{\circ r}$. From this, we deduce that the almost-union of the $a_i B_R$ constructed in the previous paragraph contains each $a_i B_{R}^{\circ rk} = a_i B_{R - rk}$. By assumption their union covers $B_{R+1}$, and by induction we can show that every $B_{r+j}$ can be constructed.

Now observe that in a general finitely-generated group, we have $B_{R+1} = \bigcup_{a \in B_{rk + 1}} a B_{R - rk}$ for all large enough $R$. So if the property fails, we must have $|B_{rk + 1}| \geq 2^k$ for infinitely many $k$. In other words, the group has exponential growth.\footnote{Recall that the size of balls is submultiplicative, so $|B_{r k + 1}| \geq 2^k$ for infinitely many $k$ implies $B_k \geq \alpha^k$ for some $\alpha > 1$ and all $k \in \N$.}
\end{proof}

Note that the free group has finite asymptotic dimension and exponential growth, while the Grigorchuk group has infinite asymptotic dimension and subexponential growth.

\section*{Acknowledgements}

The notion of homotopy arose as a tangent of the second-named author's work with Thomas Worsch in 2016. The observation that strong contractibility corresponds to equiconnectedness was given to us by user Tyrone on MathOverflow  \cite{TyMO23}. Theorem~\ref{thm:CohomologicalExample} is roughly based on a geometric idea of Ilkka T\"orm\"a on the group $\Z^2$. We thank the anonymous referee for detailed comments, which improved the paper significantly.


\printbibliography{}

\end{document}